\documentclass[12pt, reqno]{article}
\usepackage[utf8]{inputenc}
\usepackage{amsmath, amsfonts, amssymb, amsthm, mathtools, comment, float}
\usepackage[margin=1.5in]{geometry}
\usepackage[breaklinks=True]{hyperref}
\DeclareMathOperator{\id}{id}

\DeclareMathOperator{\Gal}{Gal}
\DeclareMathOperator{\tr}{tr}
\newcommand{\set}[1]{\left\{ #1 \right\}}
\newcommand{\abs}[1]{\left| #1 \right|}
\newcommand{\tup}[1]{\left( #1 \right)}
\newcommand{\ive}[1]{\left[ #1 \right]}
\newcommand{\Z}{\mathbb{Z}}

\newcommand{\Q}{\mathbb{Q}}
\newcommand{\Frob}{\mathbf{F}}
\newcommand{\G}{\mathcal{G}}
\newcommand{\one}{\underline{1}}

\newcommand{\seqnum}[1]{\href{https://oeis.org/#1}{\rm \underline{#1}}}

\newcommand{\email}[1]{\href{mailto:#1}{\texttt{#1}}}

\newtheoremstyle{plainsl}
  {8pt plus 2pt minus 4pt}
  {8pt plus 2pt minus 4pt}
  {\slshape}
  {0pt}
  {\bfseries}
  {.}
  {5pt plus 1pt minus 1pt}
  {}

\theoremstyle{plainsl}
  \newtheorem{theorem}{Theorem}[section]
  \newtheorem{proposition}[theorem]{Proposition}
  \newtheorem{lemma}[theorem]{Lemma}
  \newtheorem{corollary}[theorem]{Corollary}
\theoremstyle{definition}
  \newtheorem{definition}[theorem]{Definition}
  
\theoremstyle{remark}

\newenvironment{statement}{\begin{quote}}{\end{quote}}

\title{Necklaces over a Group with Identity Product}
\author{Darij Grinberg\footnote{Department of Mathematics, Drexel University,
(\email{darijgrinberg@gmail.com})}
\ and Peter Mao\footnote{Department of Mathematics, Drexel University,
(\email{pm929@drexel.edu})}}
\date{4 August 2026}
\begin{document}

\maketitle

\begin{abstract}
We address two variants of the classical necklace counting problem from enumerative combinatorics.
In both cases, we fix a finite group $\mathcal{G}$ and a positive integer $n$.
In the first variant, we count the ``identity-product $n$-necklaces'' --- that is, the orbits of $n$-tuples $\left(a_1, a_2, \ldots, a_n\right) \in \mathcal{G}^n$ that satisfy $a_1 a_2 \cdots a_n = 1$ under cyclic rotation.
In the second, we count the orbits of all $n$-tuples $\left(a_1, a_2, \ldots, a_n\right) \in \mathcal{G}^n$ under cyclic rotation and left multiplication (i.e., the operation of $\mathcal{G}$ on $\mathcal{G}^n$ given by $h \cdot \left(a_1, a_2, \ldots, a_n\right) = \left(ha_1, ha_2, \ldots, ha_n\right)$).
We prove bijectively that both answers are the same, and express them as a sum over divisors of $n$.

Consequently, we generalize the first problem to $n$-necklaces whose product of entries lies in a given subset of $\mathcal{G}$ (closed under conjugation), and we connect a particular case to the enumeration of irreducible polynomials over a finite field with given degree and second-highest coefficient $0$.
\medskip

\medskip

\end{abstract}

\section{Introduction}

Necklaces are combinatorial objects situated on the spectrum between ordered selections (tuples) and unordered selections (multisets) of elements from a set.
Formally, they are defined as orbits (i.e., equivalence classes) of $n$-tuples under cyclic rotation. They can be visualized as labellings (or colorings) of a circular necklace, where we understand that the necklace can be rotated at will (but not turned over) without changing the object.

Counting necklaces is a classical problem in enumerative combinatorics (see, e.g., \cite[\S 3, Example 3 (b)]{BenGol}, \cite[Chapter 7]{Loehr}, \cite[\S 5.5]{GolLiu}).
A well-known formula due to Moreau \cite{Moreau72} says that the number of $n$-necklaces (i.e., rotation classes of $n$-tuples) over a size-$q$ set $A$ is
\begin{equation}
\dfrac{1}{n} \sum_{d\mid n} \phi\left(d\right) q^{n/d},
\label{eq.neck}
\end{equation}
where the sum runs over all positive divisors $d$ of $n$, and where $\phi$ denotes the Euler totient function (see Definition~\ref{def.phi} below).
This surprisingly simple formula is a starting point for several deeper developments, including the theory of necklace polynomials \cite[\S 17.3]{Haz-Witt}.
The fact that \eqref{eq.neck} is an integer (obvious from its combinatorial meaning) itself has some uses in number theory (yielding, e.g., Fermat's Little Theorem when $n$ is prime).
The expression \eqref{eq.neck} also counts irreducible monic polynomials of degree dividing $n$ over the finite field $\mathbb{F}_q$ when $q$ is a prime power, and this can be proved using a bijection between necklaces and irreducible polynomials \cite[\S 7.6.2]{Reutenauer}.
There is also a bijection between words and multisets of necklaces \cite[\S 3]{GesReu} with far-reaching applications to symmetric functions, free Lie algebras and permutation enumeration.
Further afield, the Burrows--Wheeler transform (transforming words into necklaces) is used for data compression in 
software such as \texttt{bzip2} \cite[\S 2.3]{Adjeroh2008}.

In this paper, we count what we call \emph{identity-product necklaces}, i.e., necklaces where instead of colored beads we have elements of a finite group $\G$ which multiply to the identity.
This is an orbit-counting problem for a certain action of the cyclic group $C_n$, which (by a well-known technique) can be reduced to counting identity-product tuples which are fixed under different cyclic rotations.
We thus obtain a ``sum over divisors'' formula similar to \eqref{eq.neck} but involving the sizes of the group and of its torsion sets (Theorem~\ref{thm.main}).
Once again, the integrality of the result is not evident from the formula, and thus gives rise to nontrivial congruences between sizes of torsion sets modulo the length of the necklaces being counted.

After proving our main formula, we count the \emph{aperiodic} identity-product necklaces (i.e., the ones that have no nontrivial rotational symmetries) as well (Theorem~\ref{thm.aperiodic}).
As expected from classical necklace-counting, the sum is the same but for the replacement of the Euler totient function $\phi$ by the M\"obius function $\mu$.

The number of identity-product $n$-necklaces over a group $\G$ furthermore equals the number of ``homogeneous $n$-necklaces'' over $\G$, that is, orbits of $n$-tuples of elements of $\G$ under the action of the direct product $C_n \times \G$ (where $C_n$ acts by cyclic rotation and $\G$ acts by left multiplication on each entry).
We prove this bijectively (Theorem~\ref{thm.hg-main}), thus obtaining a formula for the latter number as well.

In Section~\ref{sect.Kprod}, we extend some of our enumerative results to the more general notion of \emph{$K$-product necklaces}, in which the product of the entries is required to belong to a given conjugacy-invariant subset $K$ of the group.

Just like Moreau's original necklaces, our identity-product necklaces can be connected to irreducible polynomials over finite fields. In Section~\ref{sect.irpol}, we take up this connection and prove some enumerative results about conjugacy classes in finite field extensions and irreducible polynomials.

In Section~\ref{sect.oeis}, we relate our enumerations to three integer sequences from the On-Line Encyclopedia of Integer Sequences.
In the short final Section~\ref{sect.open}, we suggest two venues for further research.

\section{Definitions}

Fix a finite group $\G$ with binary operation $\cdot$ and identity element $1$.

\begin{definition}
    For any nonnegative integer $n$, we let $\G_1^n$ be the set of all \emph{identity-product $n$-tuples}, defined by
    \begin{align*}
        \G_1^n &:= \left\{a = (a_1, a_2, \dots, a_n) \in \G^n : a_1 a_2 a_3 \cdots a_n = 1\right\} .
    \end{align*}
\end{definition}

\begin{definition}
    For any integer $n$, let us define the number $[\G \div n]$ to be
    \begin{align*}
        [\G \div n] & := \abs{\left\{a \in \G : a^n = 1\right\}} .
    \end{align*}
    This is the number of elements in $\G$ whose order divides $n$.
\end{definition}

\begin{definition}
    Let $n$ be a positive integer.
    Then, the $n$-element cyclic group
	\[
	C_n = \set{1, g, g^2, \ldots, g^{n-1}}
	\]
	acts on $\G^n$ by cyclic rotation:
    \[
    g\tup{a_1, a_2, \ldots, a_n} = \tup{a_n, a_1, a_2, \ldots, a_{n-1}}.
    \]
    To see that the subset $\G_1^n$ of $\G^n$ is preserved by this action, we observe that the equation
    \begin{align*}
        a_1 a_2 a_3 \cdots a_{n-1} a_n &= 1
    \end{align*}
    implies that $a_n$ is the inverse of $a_1 a_2 a_3 \cdots a_{n-1}$, and thus
    \begin{align*}
        a_n a_1 a_2 a_3 \cdots a_{n-1} &= 1 .
    \end{align*}
    Therefore, if $a = \tup{a_1, a_2, a_3, \dots, a_n}$ is an identity-product $n$-tuple, then so is $ga$. Hence $\G_1^n$ is preserved by the action of $C_n$ and thus itself becomes a $C_n$-set.
    The orbits of this action of $C_n$ on $\G_1^n$ (that is, the elements of the quotient $\G_1^n / C_n$) are called the \emph{identity-product $n$-necklaces}.
\end{definition}

\begin{definition}
\label{def.phi}
    Let $n$ be a positive integer. We define the \emph{Euler totient} $\phi\tup{n}$ to be the number of integers $i \in \set{1, 2, \ldots, n}$ that are coprime to $n$.
\end{definition}

\section{Main theorem I: Identity-product $n$-necklaces}

Our first main result is a formula for the number of identity-product $n$-necklaces for any finite group $\G$:

\begin{theorem}
    \label{thm.main}
    Let $n$ be a positive integer.
    The number $\abs{\G_1^n / C_n}$ of identity-product $n$-necklaces is
    \begin{align*}
        \abs{\G_1^n/C_n}
        &= \frac{1}{n}\sum_{d \mid n} \phi\left(\frac{n}{d}\right) \left[\G \div \frac{n}{d}\right] \abs{\G}^{d-1} .
    \end{align*}
    Here, as usual, the summation sign means a sum over all positive divisors $d$ of $n$.
\end{theorem}

We can make this more explicit if $\G$ is abelian. Recall the cyclic groups $C_k$ of order $k$, defined for all positive integers $k$.
The structure theorem for finite abelian groups says that each abelian group is isomorphic to a direct product of such cyclic groups.
If $\G$ is an abelian group written in such a form, then the above formula reduces to the following:

\begin{corollary}
    \label{cor.main}
    Assume that
    \begin{align*}
        \G &= \prod_{i=1}^m C_{k_i}
        \qquad \tup{\text{a direct product of cyclic groups}}
    \end{align*}
    for some positive integers $k_1, k_2, \ldots, k_m$.
    Then the number $\abs{\mathcal{G}_1^n / C_n}$ of identity-product $n$-necklaces is
    \begin{align*}
        \abs{\mathcal{G}_1^n/C_n} &= \frac{1}{n}\sum_{d \mid n} \phi\left(\frac{n}{d}\right) \prod_{i=1}^m \left(\gcd\left(k_i, \frac{n}{d}\right) k_i^{d-1}\right) .
    \end{align*}
\end{corollary}

\section{Lemmas}

Our proofs of the above claims are based on a common set of lemmas.
We use the \emph{orbit-counting lemma} (by historic accident also known as the \emph{Burnside lemma}; see \cite{Vatter20} or \cite[Theorem 6.2.5]{Guichard} for a proof):

\begin{lemma}
    \label{lem.burnside}
    Let $G$ be a finite group acting on a finite set $X$. Then, the number $\abs{X / G}$ of orbits is
    \begin{align*}
        \abs{X / G} &= \frac{1}{\abs{G}}\sum_{g \in G} \abs{\{x \in X : g \cdot x = x\}} .
    \end{align*}
\end{lemma}

We also need a basic theorem due to Gauss \cite[Theorem 2.40]{Shoup2}:

\begin{lemma}
    \label{lem.gauss}
    For any positive integer $n$, we have 
    \begin{align*}
        \sum_{d \mid n} \phi\tup{d} = n.
    \end{align*}
\end{lemma}

A useful consequence thereof is the following:

\begin{lemma}
    \label{lem.gaussPlusDirichlet}
    Let $f$ and $g$ be two functions from $\set{1,2,3, \dots}$ to $\Q$. Assume that
    \begin{align}
        g(k) &= \sum_{d \mid k} f(d) \qquad \text{for each } k \geq 1. & \label{eq.gaussPlusDirichletAssumption}
    \end{align}
    Then
    \begin{align}
        \sum_{d \mid k} \frac{k}{d}f(d) &= \sum_{d \mid k} \phi\left(\frac{k}{d}\right) g\left(d\right) \qquad \text{for each } k \geq 1. & \label{eq.gaussPlusDirichletConclusion}
    \end{align}
\end{lemma}

We shall prove this lemma using the following classical notion:

\begin{definition}
    Let $f$ and $g$ be two functions from $\set{1,2,3,\dots}$ to $\Q$. We define the \emph{Dirichlet convolution} $f * g$ to be the function from $\set{1,2,3,\dots}$ to $\Q$ given by
    \begin{align*}
        \tup{f * g} \tup{k} &= \sum_{d \mid k} f\tup{d} g\left(\frac{k}{d}\right) .
    \end{align*}
\end{definition}

The Dirichlet convolution operation $*$ is known to be commutative and associative (see, e.g., \cite[\S 2.9]{Shoup2}).
Furthermore, Lemma~\ref{lem.gauss} can be restated as 
\begin{align}
    \phi * \one = \id , \label{eq.phione}
\end{align}
where $\one$ denotes the constant function that sends all positive integers to $1$ (because every integer $n > 0$ satisfies $(\phi * \one)(n) = \sum_{d\mid n} \phi(d) \one(\frac{n}{d}) = \sum_{d\mid n} \phi(d)$ and $\id(n) = n$).

\begin{proof}[Proof of Lemma \ref{lem.gaussPlusDirichlet}]
    We note that the assumption \eqref{eq.gaussPlusDirichletAssumption} is the written out form of the statement
    \begin{align}
        g &= f * \one \label{eq.gequalsfone}
    \end{align}
    where $\one$ is the constant function which maps all positive integers to $1$. We can also rewrite the conclusion \eqref{eq.gaussPlusDirichletConclusion} as
    \begin{align}
        \sum_{d \mid k} f(d) \frac{k}{d} &= \sum_{d \mid k} g\left(d\right) \phi\left(\frac{k}{d}\right) ,
    \end{align}
    or, equivalently,
    \begin{align}
        (f * \id)(k) &= (g * \phi)(k) .
    \end{align}
    To prove this, we observe that
    \begin{align*}
        g * \phi &= (f * \one) * \phi & \tup{\text{by \eqref{eq.gequalsfone}}}\\
        &= f * (\one * \phi) & \tup{\text{by associativity of the Dirichlet convolution}} \\
        &= f * (\phi * \one) & \tup{\text{by commutativity of the Dirichlet convolution}} \\
        &= f * \id & \tup{\text{by \eqref{eq.phione}}}
    \end{align*}
    and therefore $f * \id = g * \phi$.
\end{proof}

The following lemma is a classical property of greatest common divisors:

\begin{lemma}
    \label{lem.gcd-div}
    Let $a$, $b$ and $c$ be three integers, where $a$ and $b$ are not both $0$. Then, $a \mid bc$ holds if and only if $\dfrac{a}{\gcd\tup{a,b}} \mid c$.
\end{lemma}

\begin{proof}
    The ``if'' part is easy (since multiplying the divisibility $\dfrac{a}{\gcd\tup{a,b}} \mid c$ with the obvious divisibility $\gcd\tup{a,b} \mid b$ yields $a \mid bc$).
    It remains to prove the ``only if'' part.
    So let us assume that $a \mid bc$, and try to prove that $\dfrac{a}{\gcd\tup{a,b}} \mid c$.

    Let $g = \gcd\left(a,b\right)$. Thus, we can write $a$ and $b$ as $a = ga'$ and $b = gb'$ for two coprime integers $a'$ and $b'$.
    Consider these $a'$ and $b'$. Thus, $a \mid bc$ rewrites as $ga' \mid gb'c$. Cancelling $g$ (which is nonzero, since $a$ and $b$ are not both $0$), we thus obtain $a' \mid b'c$.
    Since $a'$ is coprime to $b'$, we can cancel $b'$ from this divisibility (by \cite[Theorem 1.9]{Shoup2}), and obtain $a' \mid c$.
    In other words,
	\[
	\dfrac{a}{\gcd\tup{a,b}} \mid c
	\]
    (since $a' = \dfrac{a}{g} = \dfrac{a}{\gcd\tup{a,b}}$), which completes our proof.
\end{proof}

\begin{lemma}
    \label{lem.torsion-num}
    Let $a$ and $b$ be two positive integers. Then,
    $\ive{ C_a \div b } = \gcd\tup{a, b}$.
\end{lemma}

\begin{proof}
    Set $g = \gcd\tup{a,b}$ and $a' = a / g$, so that $a = ga'$.
    We write $\overline{c}$ for the residue class of an integer $c$ modulo $a$. This residue class belongs to the finite ring $\Z / a \Z$.
    
    Now, by the definition of $\ive{ C_a \div b }$, we have
    \begin{align*}
        \ive{ C_a \div b }
        &=
        \abs{\set{h \in C_a : h^b = 1}}\\
        &=
        \abs{\set{x \in \Z / a \Z : bx = \overline{0}}}
    \end{align*}
	(since $C_a$ is the additive abelian group $\Z / a \Z$, rewritten multiplicatively).
	Hence,
    \begin{align}
        \ive{ C_a \div b }
        &=
		\abs{\set{x \in \Z / a \Z : bx = \overline{0}}} \nonumber\\
        &=
        \abs{\set{c \in \set{0, 1, \ldots, a-1} : b\overline{c} = \overline{0}}}
        \label{pf.lem.torsion-num.1}
    \end{align}
    (since the elements of $\Z / a \Z$ are precisely the elements $\overline{c}$ with $c \in \set{0, 1, \ldots, a-1}$).
    However, for any $c \in \set{0, 1, \ldots, a-1}$, we have the chain of equivalences
    \begin{align*}
        \tup{b\overline{c} = \overline{0}}
        &\Longleftrightarrow
        \tup{bc \equiv 0 \mod a}
        \Longleftrightarrow
        \tup{a \mid bc} \\
        &\Longleftrightarrow
        \tup{\dfrac{a}{\gcd\tup{a,b}} \mid c}
        \qquad \tup{\text{by Lemma~\ref{lem.gcd-div}}}\\
        &\Longleftrightarrow
        \tup{\dfrac{ga'}{g} \mid c}
        \qquad \tup{\text{since $a = ga'$ and $\gcd\tup{a,b} = g$}}\\
        &\Longleftrightarrow
        \tup{a' \mid c}.
    \end{align*}
    Thus, we can rewrite \eqref{pf.lem.torsion-num.1} as
    \begin{align*}
        \ive{ C_a \div b }
        &=
        \abs{\set{c \in \set{0, 1, \ldots, a-1} : a' \mid c}} \\
        &=
        \abs{\set{c \in \set{0, 1, \ldots, ga'-1} : a' \mid c}}
        \qquad \tup{\text{since } a = ga'} \\
        &=
        \abs{\set{\text{the multiples of $a'$ between $0$ (inclusive) and $ga'$ (exclusive)}}}\\
        &=
        \abs{\set{0a', 1a', 2a', \ldots, \tup{g-1}a'}}
        = g = \gcd\tup{a,b}.
    \end{align*}
\end{proof}

\section{Proofs}

We now come to the proofs of Theorem \ref{thm.main} and Corollary \ref{cor.main}.

\begin{proof}[Proof of Theorem \ref{thm.main}]
    A \emph{period} of an $n$-tuple $a = (\alpha_1, \alpha_2, \dots, \alpha_n)$ means a positive integer $k \in \set{1, 2, 3, \dots}$ such that $g^k \cdot a = a$.
    Clearly, $n$ is a period of each $n$-tuple.
    Thus, each $n$-tuple $a$ has a smallest period.
    This smallest period must furthermore divide $n$ (because if it does not, then the remainder of $n$ upon division by this smallest period will be an even smaller period of $a$, which is a contradiction).
    
    Next, we define $P_k$ to be the set of all identity-product tuples in $\G^n_1$ which have smallest period $k$:
    \begin{align}
        P_k &:= \set{a \in \G^n_1 : g^k \cdot a = a,\text{ but }  g^l \cdot a \neq a \text{ for all } 0 < l < k } .
        \label{def.Pk}
    \end{align}
    Finally, we define $A_k$ to be the set of all identity-product tuples with period $k$:
    \begin{align}
        A_k &:= \set{a \in \G^n_1 : g^k \cdot a = a} .
    \end{align}
    
    \begin{lemma}
        \label{lem.sizeOfAi}
        For any $k \mid n$, we have
        \begin{align}
            \abs{A_k} &= \left[\G \div \frac{n}{k}\right] \abs{\G}^{k-1} .
        \end{align}
    \end{lemma}
    
    \begin{proof}[Proof of Lemma \ref{lem.sizeOfAi}]
        By definition, $A_k$ is the set of identity-product tuples fixed by $g^k$. This implies that if we denote the first $k$ entries of $a \in A_k$ by $\alpha_1, \alpha_2, \dots, \alpha_k$, then
        \begin{align*}
            a &= (\alpha_1, \alpha_2, \dots, \alpha_k, \underbrace{*, \dots, *}_{(n-k) \text{ entries}})
            \\
            &= g^k (\alpha_1, \alpha_2, \dots, \alpha_k, \underbrace{*, \dots, *}_{(n-k) \text{ entries}}) \\
            &= (\underbrace{*, \dots, *}_{k \text{ entries}}, \alpha_1, \alpha_2, \dots, \alpha_k,\underbrace{*, \dots, *}_{(n-2k) \text{ entries}})
            \\
            &= g^k(\underbrace{*, \dots, *}_{k \text{ entries}}, \alpha_1, \alpha_2, \dots, \alpha_k,\underbrace{*, \dots, *}_{(n-2k) \text{ entries}}) \\
            &\qquad \qquad \vdots \\
            &= (\underbrace{*, \dots, *}_{(n-k) \text{ entries}}, \alpha_1, \alpha_2, \dots, \alpha_k) .
        \end{align*}
        Comparing every other line entrywise, we conclude that
        \begin{align*}
            a &= (\alpha_1, \alpha_2, \dots, \alpha_k, \alpha_1, \alpha_2, \dots, \alpha_k, \dots, \alpha_1, \alpha_2, \dots, \alpha_k) .
        \end{align*}
        This means that the whole identity-product tuple $a \in A_k$ is determined by $\alpha_1, \alpha_2, \dots, \alpha_k$.
        Moreover, these $k$ elements $\alpha_1, \alpha_2, \dots, \alpha_k$ must satisfy
        \begin{align*}
            (\alpha_1 \cdot \alpha_2 \cdot \ldots \cdot \alpha_k)^{n/k} = 1
        \end{align*}
        (in order for the product of the entries of $a$ to be $1$).
        We can thus construct such a tuple $a$ by first choosing the $k-1$ entries $\alpha_1, \alpha_2, \dots, \alpha_{k-1}$ arbitrarily (there are $\abs{\G}^{k-1}$ many options for this), and then choosing the product $\alpha_1 \cdot \alpha_2 \cdot \ldots \cdot \alpha_k$ to be any element $h \in \G$ that satisfies $h^{n/k} = 1$ (there are $\left[\G \div \frac{n}{k}\right]$ many such elements); this uniquely determines $\alpha_k$ and therefore the entire tuple $a$.
        So the total number of such $n$-tuples $a \in A_k$ is $\abs{\G}^{k-1} \cdot \left[\G \div \frac{n}{k}\right]$. This proves Lemma \ref{lem.sizeOfAi}.
    \end{proof}
    
    \begin{lemma}
        \label{lem.relateAtoP}
        For any positive integer $k$, we have
        \begin{align}
            \abs{A_k} &= \sum_{d \mid k} \abs{P_d}.
        \end{align}
    \end{lemma}
    
    \begin{proof}[Proof of Lemma \ref{lem.relateAtoP}]
        Each $a \in A_k$ has some smallest period $d$,
        and this $d$ must divide $k$ (since otherwise, the remainder of $k$ upon division by $d$ would be an even smaller period of $a$).
        In other words, each $a \in A_k$ belongs to $P_d$ for some $d \mid k$, and this $d$ is obviously unique.
        Thus,
        \begin{align}
            A_k \subseteq \bigsqcup_{d \mid k} P_d .
            \label{pf.lem.relateAtoP.sub}
        \end{align}
        Conversely, if $p \in P_d$ for some $d \mid k$, then the tuple $p$ has smallest period $d$ and thus has period $k$ (since $k$ is a multiple of $d$), so that $p$ belongs to $A_k$.
        Thus,
        \[
                \bigsqcup_{d \mid k} P_d \subseteq A_k.
        \]
        Combined with \eqref{pf.lem.relateAtoP.sub}, this leads to
        \begin{align}
            A_k &= \bigsqcup_{d \mid k} P_d .
        \end{align}
        Hence,
        \[
        \abs{A_k} = \abs{\bigsqcup_{d \mid k} P_d} = \sum_{d \mid k} \abs{P_d} .
        \]
        This proves Lemma \ref{lem.relateAtoP}.
    \end{proof}
    
    Now Lemma \ref{lem.burnside} (applied to $G = C_n$ and $X = \G_1^n$) yields
    \begin{align}
        \abs{\G_1^n / C_n} &= \frac{1}{\abs{C_n}}\sum_{h \in C_n} \abs{\set{a \in \G_1^n : h \cdot a = a}} \nonumber\\
        &= \frac{1}{n} \sum_{k=1}^{n} \abs{\set{a \in \G_1^n : g^k \cdot a = a}}\nonumber
    \end{align}
    (since $C_n$ consists of the $n$ elements $g^1, g^2, \ldots, g^n$).
    Since each $k \in \set{1, 2, \dots, n}$ satisfies
    \begin{align}
        \abs{\set{a \in \G_1^n : g^k \cdot a = a}} &= \abs{A_k} = \sum_{d \mid k} \abs{P_d} & \tup{\text{by Lemma \ref{lem.relateAtoP}}},\nonumber
    \end{align}
    we can rewrite this as
    \begin{align*}
        \abs{\G_1^n / C_n} &= \frac{1}{n} \sum_{k=1}^{n} \sum_{d \mid k} \abs{P_d} .
    \end{align*}
    Since the smallest period of an $n$-tuple must divide $n$, we have $\abs{P_d} = 0$ whenever $d \nmid n$. This gives us
    \begin{align*}
        \frac{1}{n} \sum_{k=1}^{n} \sum_{d \mid k} \abs{P_d} &= \frac{1}{n} \sum_{k=1}^{n} \sum_{\substack{d \mid k; \\ d \mid n}} \abs{P_d} .
    \end{align*}
    Now we switch the order of the summations:
    \begin{align*}
        \frac{1}{n} \sum_{k=1}^{n} \sum_{\substack{d \mid k; \\ d \mid n}} \abs{P_d} &= \frac{1}{n} \sum_{d \mid n} \sum_{\substack{k \in [1,n]; \\ d \mid k}} \abs{P_d} .
    \end{align*}
    Then we note that there are $\frac{n}{d}$ equal terms in the inner sum since $d \mid k$ is equivalent to $k$ being a multiple of $d$, so
    \begin{align*}
        \frac{1}{n} \sum_{d \mid n} \sum_{\substack{k \in [1,n]; \\ d \mid k}} \abs{P_d} &= \frac{1}{n} \sum_{d \mid n} \frac{n}{d} \abs{P_d} = \sum_{d \mid n} \frac{\abs{P_d}}{d} .
    \end{align*}
    Combining the previous five equalities, we obtain
    \begin{align}
        \abs{\G_1^n / C_n} &= \sum_{d \mid n} \frac{\abs{P_d}}{d} \label{pf.lem.relateAtoP.2} .
    \end{align}
    Let us define two functions $f$ and $g$ from $\set{1, 2, 3, \ldots}$ to $\Z$ by
    \begin{align}
        f(k) &= \abs{P_k}; \\
        g(k) &= \abs{A_k} .
    \end{align}
    Then, the assumption of Lemma~\ref{lem.gaussPlusDirichlet} is satisfied because of Lemma~\ref{lem.relateAtoP}. Hence, the former lemma yields
    \[
    \sum_{d \mid n} \frac{n}{d}f(d)
        = \sum_{d \mid n} \phi\left(\frac{n}{d}\right)g(d) .
    \]
    In view of our definitions of $f$ and $g$ and dividing both sides by $n$, we can rewrite this as
    \[
    \sum_{d \mid n} \frac{\abs{P_d}}{d}
        = \frac{1}{n}\sum_{d \mid n} \phi\left(\frac{n}{d}\right)\abs{A_d} .
    \]
    Using \eqref{pf.lem.relateAtoP.2} for the left hand side and Lemma \ref{lem.sizeOfAi} for the right hand side, we can further rewrite this as
    \begin{align*}
        \abs{\G_1^n / C_n} &= \frac{1}{n}\sum_{d \mid n} \phi\left(\frac{n}{d}\right)\left[\G \div \frac{n}{d}\right] \abs{\G}^{d-1} .
    \end{align*}
    This proves Theorem \ref{thm.main}.
\end{proof}

\begin{proof}[Proof of Corollary \ref{cor.main}]
    To begin with, since
    \begin{align}
        \G &= \prod_{i=1}^m C_{k_i} ,
        \label{eq.G=prod}
    \end{align}
    we can compute the size
    \begin{align}
        \abs{\G} &= \abs{\prod_{i=1}^m C_{k_i}}
        = \prod_{i=1}^m \abs{C_{k_i}}
        = \prod_{i=1}^m k_i .
        \label{pf.cor.main.1}
    \end{align}
    Now we need to compute $\left[\G \div \frac{n}{d}\right]$ for each $d \mid n$. Any finite groups $H_1, H_2, \ldots, H_m$ and any integer $k$ satisfy
    \[
    \ive{ \prod_{i=1}^m H_i \div k } = \prod_{i=1}^m \ive{ H_i \div k},
    \]
    since $k$-th powers in a direct product are taken entrywise.
    Thus, from \eqref{eq.G=prod}, we obtain
    \begin{align}
        \left[\G \div \frac{n}{d}\right] &= \prod_{i=1}^m \left[C_{k_i} \div \frac{n}{d}\right] = \prod_{i=1}^m \gcd\tup{k_i, \frac{n}{d}}
        \label{pf.cor.main.2}
    \end{align}
    (where we used Lemma \ref{lem.torsion-num} in the last step).
    The claim of Theorem \ref{thm.main} thus becomes
    \begin{align*}
        \abs{\G^n_1 / C_n} &= \frac{1}{n} \sum_{d \mid n} \phi\left(\frac{n}{d}\right) \ive{\G \div \frac{n}{d}} \abs{\G}^{d-1} \\
        &= \frac{1}{n} \sum_{d \mid n} \phi\left(\frac{n}{d}\right) \left(\prod_{i=1}^m \gcd\tup{k_i, \frac{n}{d}}\right) \left(\prod_{i=1}^m k_i\right)^{d-1}
        & \tup{\text{by \eqref{pf.cor.main.1} and \eqref{pf.cor.main.2}}}
        \\
        &= \frac{1}{n} \sum_{d \mid n} \phi\left(\frac{n}{d}\right) \prod_{i=1}^m \left(\gcd\tup{k_i, \frac{n}{d}} k_i^{d-1}\right) .
    \end{align*}
\end{proof}

\section{\label{sec.aperiodic}Main theorem II: Aperiodic necklaces}

After having computed the total number of identity-product $n$-necklaces for a given finite group $\G$, we can further restrict ourselves to just those necklaces which are aperiodic.

We declare an identity-product $n$-necklace $N$ to be \emph{aperiodic} if it contains $n$ distinct $n$-tuples. In other words, a necklace $N = \set{a, g\cdot a, g^2\cdot a, \ldots, g^{n-1}\cdot a} \subseteq \G_1^n$ is aperiodic if and only if no two of the tuples $a, g\cdot a, g^2\cdot a, \ldots, g^{n-1} \cdot a$ are equal. This latter condition is equivalent to $a \in P_n$, where we use the notation $P_k$ introduced in \eqref{def.Pk} (since $g^n \cdot a = a$ holds automatically).
Thus, each aperiodic necklace in $\G_1^n$ contains $n$ distinct tuples in $P_n$, and each tuple in $P_n$ is contained in a unique aperiodic necklace in $\G_1^n$.
Hence, the number of aperiodic necklaces in $\G_1^n$ is $\dfrac{\abs{P_n}}{n}$.

Therefore, it suffices to compute $\abs{P_n}$.
To compute it, we need another number-theoretic function instead of the Euler totient function:

\begin{definition}
     We define the M\"obius function $\mu : \set{1,2,3,\ldots} \to \Z$ by
    \begin{align*}
        \mu(n) &=
		\begin{cases}
		(-1)^m, & \text{if }n \text{ is a product of }m\text{ distinct primes;} \\
		0, & \text{otherwise.}
		\end{cases}
    \end{align*}
\end{definition}

\begin{theorem}
    \label{thm.aperiodic}
    Let $n$ be a positive integer. The number $\frac{\abs{P_n}}{n}$ of aperiodic identity-product $n$-necklaces is
    \begin{align*}
        \frac{\abs{P_n}}{n} &= \frac{1}{n}\sum_{d \mid n} \mu\left(\frac{n}{d}\right)\ive{\G \div \frac{n}{d}}\abs{\G}^{d-1} .
    \end{align*}
\end{theorem}

Our proof of this theorem relies on the \emph{M\"obius inversion formula}, a classical result in elementary number theory (see, e.g., \cite[Theorem 2.39 and the paragraph that follows]{Shoup2}):
\begin{lemma}
    \label{lem.mobius}
    Let $f$ and $g$ be two functions from $\set{1, 2, 3, \ldots}$ to $\Q$.
    Assume that
    \begin{align*}
        g(k) &= \sum_{d \mid k} f(d) \qquad \text{for each } k \geq 1.
    \end{align*}
    Then 
    \begin{align*}
        f(k) &= \sum_{d \mid k} \mu(d) g\left(\frac{k}{d}\right) \qquad \text{for each } k \geq 1.
    \end{align*}
\end{lemma}
\begin{proof}[Proof of Theorem \ref{thm.aperiodic}]
Let us define two functions $f$ and $g$ from $\set{1,2,3,\dots}$ to $\Z$ by
\begin{align}
    f(k) &= \abs{P_k}, \\
    g(k) &= \abs{A_k}.
\end{align}
Then, the assumption of Lemma \ref{lem.mobius} is satisfied because of Lemma \ref{lem.relateAtoP}. Hence, Lemma \ref{lem.mobius} yields
\begin{align}
    \abs{P_k} &= \sum_{d \mid k} \mu(d) \abs{A_{\frac{k}{d}}}
    \label{eq.aperiodic.general}
\end{align}
for all $k \geq 1$.
Dividing both sides by $n$ and setting $k = n$, this becomes
\begin{align*}
    \frac{\abs{P_n}}{n} &= \frac{1}{n}\sum_{d \mid n} \mu(d) \abs{A_{\frac{n}{d}}} .
\end{align*}
We can substitute $d$ for $\frac{n}{d}$ in the sum on the right, since the map $d \mapsto \frac{n}{d}$ is a bijection from the set of positive divisors of $n$ to itself. This yields
\begin{align}
    \frac{\abs{P_n}}{n} &= \frac{1}{n}\sum_{d \mid n} \mu\left(\frac{n}{d}\right) \abs{A_{d}} .
\end{align}
Finally, we use Lemma \ref{lem.sizeOfAi} to rewrite this as
\begin{align}
    \frac{\abs{P_n}}{n} &= \frac{1}{n}\sum_{d \mid n} \mu\left(\frac{n}{d}\right) \ive{\G \div \frac{n}{d}} \abs{\G}^{d-1}
\end{align}
which concludes the proof.
\end{proof}

On top of the special case of aperiodic necklaces, let us consider the case of identity-product $n$-necklaces with smallest period $k$.
Here we declare an $n$-necklace $N$ to have \emph{smallest period} $k$ when $\abs{N} = k$.
Equivalently, an $n$-necklace $N = \{a, g\cdot a, g^2 \cdot a, \dots, g^{n-1} \cdot a\} \subseteq \G_1^n$ has smallest period $k$ if and only if the $n$-tuple $a$ has smallest period $k$.
This is because if $k$ is the smallest period of an $n$-tuple $a$, then the first $k$ of the $n$-tuples $a, g \cdot a, g^2 \cdot a, \ldots, g^{n-1} \cdot a$ are distinct, while all following tuples are repetitions of them.
Thus, each identity-product $n$-necklace with smallest period $k$ contains $k$ distinct tuples in $P_k$.
Hence, the number of these necklaces is $\frac{\abs{P_k}}{k}$.

\begin{corollary}
    \label{cor.aperiodic}
    Let $n$ and $k$ be positive integers with $k \mid n$. The number $\frac{\abs{P_k}}{k}$ of identity-product $n$-necklaces with smallest period $k$ is
    \begin{align*}
        \frac{\abs{P_k}}{k} &= \frac{1}{k}\sum_{d \mid k} \mu\left(\frac{k}{d}\right)\ive{\G \div \frac{n}{d}}\abs{\G}^{d-1} .
    \end{align*}
\end{corollary}

\begin{proof}[Proof of Corollary \ref{cor.aperiodic}]
Recall that in the proof of Theorem \ref{thm.aperiodic}, we proved the general equation \eqref{eq.aperiodic.general}. That is,
\begin{align*}
    \abs{P_k} &= \sum_{d \mid k} \mu(d) \abs{A_{\frac{k}{d}}}.
\end{align*}
By bijective substitution of $d$ with $\frac{k}{d}$ on the positive divisors of $k$ to themselves, this can be turned into
\begin{align*}
    \abs{P_k} &= \sum_{d \mid k} \mu\left(\frac{k}{d}\right) \abs{A_{d}} .
\end{align*}
We can rewrite the right hand side using Lemma \ref{lem.sizeOfAi}, applied to $d$ instead of $k$ (since $d \mid k \mid n$). The above equality thus takes the form
\begin{align}
    \abs{P_k} &= \sum_{d \mid k} \mu\left(\frac{k}{d}\right) \ive{\G \div \frac{n}{d}} \abs{\G}^{d-1} .
\end{align}
Now, dividing both sides by $k$, we obtain
\begin{align}
    \frac{\abs{P_k}}{k} &= \frac{1}{k}\sum_{d \mid k} \mu\left(\frac{k}{d}\right) \ive{\G \div \frac{n}{d}} \abs{\G}^{d-1} .
\end{align}
\end{proof}

\section{Main theorem III: Homogeneous necklaces}

An interesting variant of necklaces was proposed by Gilbert and Riordan
\cite{GilRio61}. They count $n$-tuples of elements of the cyclic group $C_{q}$
(which they write additively, as $\mathbb{Z}/q$), but identifying each
$n$-tuple not only with its cyclic rotations, but also with all the other
$n$-tuples that can be obtained by multiplying all its entries by the same
element of $C_{q}$ (in our multiplicative notation).
For instance, for $q=3$ and $n=3$, each $n$-tuple $\left(
a,b,c\right)  \in C_{q}^{n}=C_{3}^{3}$ thus ends up in a single equivalence
class with $\left(  b,c,a\right)  $, $\left(  c,a,b\right)  $, $\left(  ga,
gb, gc \right)  $, $\left(  gb, gc, ga \right)  $, $\left(  gc, ga, gb
\right)  $, $\left(  g^{2}a, g^{2}b, g^{2}c \right)  $, $\left(  g^{2}b,
g^{2}c, g^{2}a \right)  $ and $\left(  g^{2}c, g^{2}a, g^{2}b \right)  $. Such
equivalence classes can be called \textquotedblleft homogeneous
necklaces\textquotedblright. They were first considered for $q=2$ by Fine
\cite{Fine58}, motivated by psychometric experiments.

The \textquotedblleft homogeneous necklaces\textquotedblright\ can also be
generalized to an arbitrary finite group $\mathcal{G}$ instead of $C_{q}$.
Rather than introduce custom terminology for them, we shall define them as the
orbits of a group action:

As before, $\mathcal{G}$ is a finite group. We also fix a positive integer
$n$. Consider the cyclic group $C_{n}=\left\{  g^{0},g^{1},\ldots
,g^{n-1}\right\}  $ with generator $g$.

We recall that the cyclic group $C_{n}$ acts (from the left) on $\mathcal{G}^{n}$ by rotating $n$-tuples:
If $i\in\mathbb{Z}$, then the element $g^{i}\in
C_{n}$ acts on an $n$-tuple by rotating this $n$-tuple $i$ steps to the right.

The group $\mathcal{G}$ also acts from the left on $\mathcal{G}^{n}$, by the
rule
\begin{align*}
&  h\cdot\left(  a_{1},a_{2},\ldots,a_{n}\right)
= \left(  ha_{1},\ ha_{2},\ \ldots,\ ha_{n}\right) \\
&  \qquad\qquad\text{for all }h\in\mathcal{G}\text{ and }
\left(  a_{1},a_{2},\ldots,a_{n}\right)  \in\mathcal{G}^{n}.
\end{align*}
Thus, an element $h\in\mathcal{G}$ acts on an $n$-tuple $a\in\mathcal{G}^{n}$
by multiplying all entries of $a$ by $h$ from the left.

The two actions (of $C_{n}$ and of $\mathcal{G}$) on $\mathcal{G}^{n}$
commute, meaning that
\[
g^{i}\cdot\left(  h\cdot a\right)  =h\cdot\left(  g^{i}\cdot a\right)
\qquad\text{for all }g^{i}\in C_{n}\text{ and }h\in\mathcal{G}\text{ and } a\in\mathcal{G}^{n}.
\]
(This is intuitively obvious: Multiplying all entries of $a$ by $h$ can be
done either before or after a cyclic rotation; the outcome will be the same.)

This commutativity allows us to combine the actions of $C_{n}$ and of
$\mathcal{G}$ into a single action of the direct product $C_{n}\times
\mathcal{G}$ on $\mathcal{G}^{n}$. Explicitly, the group $C_{n}\times
\mathcal{G}$ acts from the left on $\mathcal{G}^{n}$ by the rule
\begin{align*}
&  \left(  g^{i},h\right)  \cdot a=g^{i}\cdot\left(  h\cdot a\right)
=h\cdot\left(  g^{i}\cdot a\right) \\
&  \qquad\qquad\text{for all }g^{i}\in C_{n}\text{ and }h\in\mathcal{G}\text{
and }a\in\mathcal{G}^{n}.
\end{align*}
Thus, a pair $\left(  g^{i},h\right)  $ acts on an $n$-tuple $a\in
\mathcal{G}^{n}$ by rotating the $n$-tuple $i$ steps to the right and
multiplying each entry by $h$ from the left. The orbits of this action can be
called the \textquotedblleft homogeneous $n$-necklaces\textquotedblright,
though we shall just refer to them as orbits. Now, we claim the following:

\begin{theorem}
\label{thm.hg-main}
The number $\abs{ \mathcal{G}^{n}/\left(  C_{n} \times\mathcal{G}\right) }$
of all orbits of the action of
$C_{n}\times\mathcal{G}$ on $\mathcal{G}^{n}$ equals the number $\left\vert
\mathcal{G}_{1}^{n}/C_{n}\right\vert $ of all identity-product $n$-necklaces.
Thus,
\[
\left\vert \mathcal{G}^{n}/\left(  C_{n}\times\mathcal{G}\right)  \right\vert
=\left\vert \mathcal{G}_{1}^{n}/C_{n}\right\vert =\frac{1}{n}\sum_{d\mid
n}\phi\left(  \frac{n}{d}\right)  \left[  \mathcal{G}\div\frac{n}{d}\right]
\abs{ \mathcal{G} }  ^{d-1}.
\]

\end{theorem}

This is implicitly claimed in the OEIS (On-Line Encyclopedia of Integer Sequences)
entry for sequence \seqnum{A000013} for the case
$\mathcal{G}=C_{2}$ (see also \newline\url{https://math.stackexchange.com/questions/3249602} ), but the general case appears new. We prove the general
case bijectively:

\begin{proof}[Proof of Theorem \ref{thm.hg-main}.]
If $i$ is any integer, then the residue
class of $i$ modulo $n$ will be denoted by $\overline{i}$. For instance,
$\overline{n} = \overline{0}$.

If $a\in\mathcal{G}^{n}$ is an $n$-tuple, then the $i$-th entry of $a$ (for
any given $i\in\left\{  1,2,\ldots,n\right\}  $) will be denoted by
$a_{\overline{i}}$ (so that $a=\left(  a_{\overline{1}},a_{\overline{2}},\ldots,a_{\overline{n}}\right)  $ and $a_{\overline{0}}=a_{\overline{n}}$
and $a_{\overline{n+1}}=a_{\overline{1}}$). Thus, the entries of our
$n$-tuples are now indexed by elements of $\mathbb{Z}/n$. If $a\in
\mathcal{G}^{n}$ is an $n$-tuple, then $a_{r}$ is hence defined for any
residue class $r\in\mathbb{Z}/n$.

Note that the action of $C_{n}$ on $\mathcal{G}^{n}$ can thus be described by
the equality
\begin{equation}
\left(  g^{i}\cdot a\right)  _{\overline{j}}=a_{\overline{j-i}}
\label{pf.thm.hg-main.actC2}%
\end{equation}
for all $g^{i}\in C_{n}$ (that is, all $i\in\mathbb{Z}$), all $a\in
\mathcal{G}^{n}$ and all $j\in\mathbb{Z}$.

Likewise, the action of $\mathcal{G}$ on $\mathcal{G}^{n}$ can be described by
the equality%
\begin{equation}
\left(  h\cdot a\right)  _{\overline{j}}=ha_{\overline{j}}
\label{pf.thm.hg-main.actG2}%
\end{equation}
for all $h\in\mathcal{G}$, all $a\in\mathcal{G}^{n}$ and all $j\in\mathbb{Z}$.

Hence, for all $g^{i}\in C_{n}$, all $h\in\mathcal{G}$, all $a\in
\mathcal{G}^{n}$ and all $j\in\mathbb{Z}$, we have%
\begin{align}
\left(  \left(  g^{i},h\right)  \cdot a\right)  _{\overline{j}}  &  =\left(
g^{i}\cdot\left(  h\cdot a\right)  \right)  _{\overline{j}}\qquad\left(
\text{since }\left(  g^{i},h\right)  \cdot a=g^{i}\cdot\left(  h\cdot
a\right)  \right) \nonumber\\
&  =\left(  h\cdot a\right)  _{\overline{j-i}}\qquad\left(  \text{by
(\ref{pf.thm.hg-main.actC2}), applied to }h\cdot a\text{ instead of }a\right)
\nonumber\\
&  =ha_{\overline{j-i}} \label{pf.thm.hg-main.actCG2}%
\end{align}
(by (\ref{pf.thm.hg-main.actG2}), applied to $j-i$ instead of $j$).

If $a\in\mathcal{G}^{n}$ is any $n$-tuple, then
\[
\left(  a_{\overline{0}}^{-1}a_{\overline{1}},\ a_{\overline{1}}%
^{-1}a_{\overline{2}},\ a_{\overline{2}}^{-1}a_{\overline{3}},\ \ldots
,\ a_{\overline{n-1}}^{-1}a_{\overline{n}}\right)  \in\mathcal{G}_{1}^{n},
\]
since
\begin{align*}
a_{\overline{0}}^{-1}a_{\overline{1}}\cdot a_{\overline{1}}^{-1}%
a_{\overline{2}}\cdot a_{\overline{2}}^{-1}a_{\overline{3}}\cdot\cdots\cdot
a_{\overline{n-1}}^{-1}a_{\overline{n}}  &  =a_{\overline{0}}^{-1}%
a_{\overline{n}}\qquad\left(  \text{by the telescope principle}\right) \\
&  =1\qquad\left(  \text{since }a_{\overline{n}}=a_{\overline{0}}\text{ (since
}\overline{n}=\overline{0}\text{)}\right)  .
\end{align*}
Thus, we can define a map%
\begin{align*}
\Delta:\mathcal{G}^{n}  &  \rightarrow\mathcal{G}_{1}^{n},\\
a  &  \mapsto\left(  a_{\overline{0}}^{-1}a_{\overline{1}},\ a_{\overline{1}%
}^{-1}a_{\overline{2}},\ a_{\overline{2}}^{-1}a_{\overline{3}},\ \ldots
,\ a_{\overline{n-1}}^{-1}a_{\overline{n}}\right)  .
\end{align*}
This map $\Delta$ satisfies the equality%
\begin{equation}
\left(  \Delta\left(  a\right)  \right)  _{\overline{j}}=a_{\overline{j-1}%
}^{-1}a_{\overline{j}} \label{pf.thm.hg-main.Del2}%
\end{equation}
for all $a\in\mathcal{G}^{n}$ and $j\in\mathbb{Z}$.

\begin{proof}[Proof of (\ref{pf.thm.hg-main.Del2}).]
Let $a\in\mathcal{G}^{n}$. The
definition of $\Delta$ shows that
\[
\Delta\left(  a\right)  =\left(
a_{\overline{0}}^{-1}a_{\overline{1}},\ a_{\overline{1}}^{-1}a_{\overline{2}%
},\ a_{\overline{2}}^{-1}a_{\overline{3}},\ \ldots,\ a_{\overline{n-1}}%
^{-1}a_{\overline{n}}\right)  .
\]
Hence, the equality $\left(  \Delta\left(
a\right)  \right)  _{\overline{j}}=a_{\overline{j-1}}^{-1}a_{\overline{j}}$
holds for each $j\in\left\{  1,2,\ldots,n\right\}  $. But both sides of this
equality --- considered as functions of $j$ --- are periodic with a period of
$n$ (since $\overline{j}$ and $\overline{j-1}$ are periodic with a period of
$n$). Therefore, having shown that this equality holds for all $j\in\left\{
1,2,\ldots,n\right\}  $, we automatically conclude that it must hold for all
$j\in\mathbb{Z}$. This proves (\ref{pf.thm.hg-main.Del2}).
\end{proof}

Next we claim the following property of $\Delta$:

\begin{statement}
\textit{Claim 1:} If two $n$-tuples $a\in\mathcal{G}^{n}$ and $b\in
\mathcal{G}^{n}$ belong to the same $C_{n}\times\mathcal{G}$-orbit, then their
images $\Delta\left(  a\right)  $ and $\Delta\left(  b\right)  $ belong to the
same $C_{n}$-orbit.
\end{statement}

\begin{proof}
[Proof of Claim 1.]Let $a\in\mathcal{G}^{n}$ and $b\in\mathcal{G}^{n}$ belong
to the same $C_{n}\times\mathcal{G}$-orbit. Thus, $b=\left(  g^{i},h\right)
\cdot a$ for some $g^{i}\in C_{n}$ and $h\in\mathcal{G}$. Consider these
$g^{i}$ and $h$. Thus, for each $j\in\mathbb{Z}$, we have%
\[
b_{\overline{j}}=\left(  \left(  g^{i},h\right)  \cdot a\right)
_{\overline{j}}=ha_{\overline{j-i}}\qquad\left(  \text{by
(\ref{pf.thm.hg-main.actCG2})}\right)
\]
and%
\[
b_{\overline{j-1}}=ha_{\overline{j-1-i}}\qquad\left(  \text{similarly}\right)
\]
and
\begin{align}
\left(  \Delta\left(  b\right)  \right)  _{\overline{j}}  &  =b_{\overline
{j-1}}^{-1}b_{\overline{j}}\qquad\left(  \text{by (\ref{pf.thm.hg-main.Del2}),
applied to }b\text{ instead of }a\right) \nonumber\\
&  =\left(  ha_{\overline{j-1-i}}\right)  ^{-1}\left(  ha_{\overline{j-i}%
}\right)  \qquad\left(  \text{since }b_{\overline{j}}=ha_{\overline{j-i}%
}\text{ and }b_{\overline{j-1}}=ha_{\overline{j-1-i}}\right) \nonumber\\
&  =\underbrace{a_{\overline{j-1-i}}^{-1}}_{=a_{\overline{j-i-1}}^{-1}%
}\underbrace{h^{-1}h}_{=1}a_{\overline{j-i}}=a_{\overline{j-i-1}}%
^{-1}a_{\overline{j-i}}. \label{pf.thm.hg-main.c1.pf.2}%
\end{align}
However, for each $j\in\mathbb{Z}$, we have%
\begin{align*}
\left(  g^{i}\cdot\Delta\left(  a\right)  \right)  _{\overline{j}}  &
=\left(  \Delta\left(  a\right)  \right)  _{\overline{j-i}}\qquad\left(
\text{by (\ref{pf.thm.hg-main.actC2}), applied to }\Delta\left(  a\right)
\text{ instead of }a\right) \\
&  =a_{\overline{j-i-1}}^{-1}a_{\overline{j-i}}\qquad\left(  \text{by
(\ref{pf.thm.hg-main.Del2}), applied to }j-i\text{ instead of }j\right)  .
\end{align*}
Comparing this with (\ref{pf.thm.hg-main.c1.pf.2}), we obtain%
\[
\left(  \Delta\left(  b\right)  \right)  _{\overline{j}}=\left(  g^{i}%
\cdot\Delta\left(  a\right)  \right)  _{\overline{j}}\qquad\text{for each
}j\in\mathbb{Z}.
\]
In other words, each entry of $\Delta\left(  b\right)  $ equals the
corresponding entry of $g^{i}\cdot\Delta\left(  a\right)  $. Hence,
$\Delta\left(  b\right)  =g^{i}\cdot\Delta\left(  a\right)  $. This shows that
$\Delta\left(  a\right)  $ and $\Delta\left(  b\right)  $ belong to the same
$C_{n}$-orbit. Claim 1 is thus proved.
\end{proof}

Next, we define a map%
\begin{align*}
\Gamma:\mathcal{G}^{n}  &  \rightarrow\mathcal{G}^{n},\\
a  &  \mapsto\left(  a_{\overline{1}},\ a_{\overline{1}}a_{\overline{2}%
},\ a_{\overline{1}}a_{\overline{2}}a_{\overline{3}},\ \ldots,\ a_{\overline
{1}}a_{\overline{2}}\cdots a_{\overline{n}}\right)  .
\end{align*}
Thus, this map $\Gamma$ is given by the rule%
\begin{equation}
\left(  \Gamma\left(  a\right)  \right)  _{\overline{j}}=a_{\overline{1}%
}a_{\overline{2}}\cdots a_{\overline{j}} \label{pf.thm.hg-main.Gam2}%
\end{equation}
for all $a\in\mathcal{G}^{n}$ and $j\in\left\{  1,2,\ldots,n\right\}  $. When
$a$ is an identity-product $n$-tuple, we can extend this rule a bit further:
For any $a\in\mathcal{G}_{1}^{n}$, we have%
\begin{equation}
\left(  \Gamma\left(  a\right)  \right)  _{\overline{j}}=a_{\overline{1}%
}a_{\overline{2}}\cdots a_{\overline{j}}\qquad\text{for all }j\in\left\{
0,1,\ldots,n\right\}  . \label{pf.thm.hg-main.Gam2a}%
\end{equation}

\begin{proof}
[Proof of (\ref{pf.thm.hg-main.Gam2a}).]Let $a\in\mathcal{G}_{1}^{n}$. Let
$j\in\left\{  0,1,\ldots,n\right\}  $. We must prove
(\ref{pf.thm.hg-main.Gam2a}). If $j\neq0$, then $j\in\left\{  1,2,\ldots
,n\right\}  $, and thus (\ref{pf.thm.hg-main.Gam2a}) follows directly from
(\ref{pf.thm.hg-main.Gam2}) in this case. Thus, it remains to prove
(\ref{pf.thm.hg-main.Gam2a}) for $j=0$. In other words, it remains to prove
that $\left(  \Gamma\left(  a\right)  \right)  _{\overline{0}}=a_{\overline
{1}}a_{\overline{2}}\cdots a_{\overline{0}}$. Since the empty product
$a_{\overline{1}}a_{\overline{2}}\cdots a_{\overline{0}}$ is $1$ by
definition, this means proving that $\left(  \Gamma\left(  a\right)  \right)
_{\overline{0}}=1$. But this is easy: We have $\overline{0}=\overline{n}$, so
that $\left(  \Gamma\left(  a\right)  \right)  _{\overline{0}}=\left(
\Gamma\left(  a\right)  \right)  _{\overline{n}}=a_{\overline{1}}%
a_{\overline{2}}\cdots a_{\overline{n}}$ (by (\ref{pf.thm.hg-main.Gam2})). But
$\left(  a_{\overline{1}},a_{\overline{2}},\ldots,a_{\overline{n}}\right)
=a\in\mathcal{G}_{1}^{n}$, and thus $a_{\overline{1}}a_{\overline{2}}\cdots
a_{\overline{n}}=1$ (by the definition of $\mathcal{G}_{1}^{n}$). Hence,
$\left(  \Gamma\left(  a\right)  \right)  _{\overline{0}}=a_{\overline{1}%
}a_{\overline{2}}\cdots a_{\overline{n}}=1$, just as we desired to prove. Thus
the proof of (\ref{pf.thm.hg-main.Gam2a}) is complete.
\end{proof}

Next, let us prove an analogue of Claim 1 for the map $\Gamma$:

\begin{statement}
\textit{Claim 2:} If two $n$-tuples $a\in\mathcal{G}_{1}^{n}$ and
$b\in\mathcal{G}_{1}^{n}$ belong to the same $C_{n}$-orbit, then their images
$\Gamma\left(  a\right)  $ and $\Gamma\left(  b\right)  $ belong to the same
$C_{n}\times\mathcal{G}$-orbit.
\end{statement}

\begin{proof}
[Proof of Claim 2.]We define a binary relation $\sim$ on the set
$\mathcal{G}^{n}$ by%
\[
\left(  u\sim v\right)  \ \Longleftrightarrow\ \left(  u\text{ and }v\text{
belong to the same }C_{n}\times\mathcal{G}\text{-orbit}\right)  .
\]
Clearly, this relation $\sim$ is an equivalence relation.

Now let us show that every $a\in\mathcal{G}_{1}^{n}$ satisfies%
\begin{equation}
\Gamma\left(  a\right)  \sim\Gamma\left(  g\cdot a\right)  .
\label{pf.thm.hg-main.Gam2.pf.1}%
\end{equation}

Indeed, fix $a\in\mathcal{G}_{1}^{n}$. We shall show that $\Gamma\left(
g\cdot a\right)  =\left(  g^{1},a_{\overline{n}}\right)  \cdot\Gamma\left(
a\right)  $ (where the pair $\left(  g^{1},a_{\overline{n}}\right)  \in
C_{n}\times\mathcal{G}$ acts on $\Gamma\left(  a\right)  $ via the
$C_{n}\times\mathcal{G}$-action on $\mathcal{G}^{n}$ defined above).

First we observe that $g\cdot a=\left(  a_{\overline{n}},a_{\overline{1}%
},a_{\overline{2}},\ldots,a_{\overline{n-1}}\right)  $ by the definition of
the $C_{n}$-action on $\mathcal{G}^{n}$. Thus, $\left(  g\cdot a\right)
_{\overline{1}}=a_{\overline{n}}$ and
\begin{equation}
\left(  g\cdot a\right)  _{\overline{i}}=a_{\overline{i-1}}\qquad\text{for
each }i\in\mathbb{Z}. \label{pf.thm.hg-main.Gam2.pf.2}%
\end{equation}

Now, let $k\in\left\{  1,2,\ldots,n\right\}  $. Then,
(\ref{pf.thm.hg-main.Gam2}) (applied to $g\cdot a$ and $k$
instead of $a$ and $j$) shows that
\begin{align}
\left(  \Gamma\left(  g\cdot a\right)  \right)  _{\overline{k}}  &  =\left(
g\cdot a\right)  _{\overline{1}}\left(  g\cdot a\right)  _{\overline{2}}%
\cdots\left(  g\cdot a\right)  _{\overline{k}} \nonumber\\
&  =\underbrace{\left(  g\cdot a\right)  _{\overline{1}}}_{=a_{\overline{n}}%
}\cdot\underbrace{\left(  g\cdot a\right)  _{\overline{2}}}%
_{\substack{=a_{\overline{1}}\\\text{(by (\ref{pf.thm.hg-main.Gam2.pf.2}))}%
}}\underbrace{\left(  g\cdot a\right)  _{\overline{3}}}%
_{\substack{=a_{\overline{2}}\\\text{(by (\ref{pf.thm.hg-main.Gam2.pf.2}))}%
}}\cdots\underbrace{\left(  g\cdot a\right)  _{\overline{k}}}%
_{\substack{=a_{\overline{k-1}}\\\text{(by (\ref{pf.thm.hg-main.Gam2.pf.2}))}%
}} \nonumber\\
&  =a_{\overline{n}}\cdot a_{\overline{1}}a_{\overline{2}}\cdots
a_{\overline{k-1}}.
\label{pf.thm.hg-main.Gam2.pf.53}
\end{align}
On the other hand,
from $k\in\left\{  1,2,\ldots,n\right\}  $,
we obtain $k-1\in\left\{  0,1,\ldots,n\right\}$.
Hence, (\ref{pf.thm.hg-main.Gam2a}) (applied to $j=k-1$) yields
$\left(  \Gamma\left(  a\right)  \right)  _{\overline{k-1}}
=a_{\overline{1}}a_{\overline{2}}\cdots a_{\overline{k-1}}$.
But
(\ref{pf.thm.hg-main.actCG2}) (applied to $1$, $a_{\overline{n}}$,
$\Gamma\left(  a\right)$ and $k$ instead of $i$, $h$, $a$ and $j$)
yields
\begin{align*}
 \left(  \left(  g^{1},a_{\overline{n}}\right)  \cdot\Gamma\left(  a\right)
\right)  _{\overline{k}}
 =a_{\overline{n}}\cdot
 \underbrace{\left(  \Gamma\left(  a\right)  \right)_{\overline{k-1}}
 }_{= a_{\overline{1}}a_{\overline{2}}\cdots a_{\overline{k-1}}}
 =a_{\overline{n}}\cdot a_{\overline{1}}a_{\overline{2}}\cdots
a_{\overline{k-1}}.
\end{align*}
Comparing this with \eqref{pf.thm.hg-main.Gam2.pf.53}, we see that
\[
\left(  \Gamma\left(  g\cdot a\right)  \right)  _{\overline{k}}=\left(
\left(  g^{1},a_{\overline{n}}\right)  \cdot\Gamma\left(  a\right)  \right)
_{\overline{k}} .
\]

Forget that we fixed $k$. We thus have shown that
$\left(  \Gamma\left(  g\cdot a\right)  \right)  _{\overline{k}}=\left(
\left(  g^{1},a_{\overline{n}}\right)  \cdot\Gamma\left(  a\right)  \right)
_{\overline{k}}$
for each $k\in\left\{  1,2,\ldots,n\right\}$.
In other words, each entry of the $n$-tuple $\Gamma\left(  g\cdot a\right)  $
equals the corresponding entry of $\left(  g^{1},a_{\overline{n}}\right)
\cdot\Gamma\left(  a\right)  $. Hence,
\[
\Gamma\left(  g\cdot a\right)  =\left(  g^{1},a_{\overline{n}}\right)
\cdot\Gamma\left(  a\right)  .
\]
Thus, the two $n$-tuples $\Gamma\left(  a\right)  $ and $\Gamma\left(  g\cdot
a\right)  $ belong to the same $C_{n}\times\mathcal{G}$-orbit. By the
definition of the relation $\sim$, this means that $\Gamma\left(  a\right)
\sim\Gamma\left(  g\cdot a\right)  $. Thus, (\ref{pf.thm.hg-main.Gam2.pf.1})
is proved.

Now, let $a\in\mathcal{G}_{1}^{n}$ and $b\in\mathcal{G}_{1}^{n}$ be two
$n$-tuples that belong to the same $C_{n}$-orbit. We must prove that their
images $\Gamma\left(  a\right)  $ and $\Gamma\left(  b\right)  $ belong to the
same $C_{n}\times\mathcal{G}$-orbit. In other words, we must prove that
$\Gamma\left(  a\right)  \sim\Gamma\left(  b\right)  $.

We have assumed that $a$ and $b$ belong to the same $C_{n}$-orbit. In other
words, $a=\gamma\cdot b$ for some $\gamma\in C_{n}$. Consider this $\gamma$.
Thus, $\gamma\in C_{n}=\left\{  g^{0},g^{1},\ldots,g^{n-1}\right\}  $. In
other words, $\gamma=g^{i}$ for some $i\in\left\{  0,1,\ldots,n-1\right\}  $.
Consider this $i$. Now, for each $k\in\mathbb{N}$, we have%
\begin{equation}
\Gamma\left(  g^{k}\cdot b\right)  \sim\Gamma\left(  g^{k+1}\cdot b\right)
\label{pf.thm.hg-main.Gam2.pf.5}%
\end{equation}
(since (\ref{pf.thm.hg-main.Gam2.pf.1}) (applied to $g^{k}\cdot b$ instead of
$a$) yields $\Gamma\left(  g^{k}\cdot b\right)  \sim\Gamma\left(
\underbrace{g\cdot g^{k}}_{=g^{k+1}}\cdot\,b\right)  =\Gamma\left(
g^{k+1}\cdot b\right)  $). Hence, we obtain a chain of relations%
\[
\Gamma\left(  g^{0}\cdot b\right)  \sim\Gamma\left(  g^{1}\cdot b\right)
\sim\Gamma\left(  g^{2}\cdot b\right)  \sim\cdots.
\]
Since $\sim$ is an equivalence relation, this yields $\Gamma\left(  g^{i}\cdot
b\right)  \sim\Gamma\left(  g^{0}\cdot b\right)  $ in particular. In view of
$g^{i}=\gamma$ and $g^{0}=1$, we can rewrite this as $\Gamma\left(
\gamma\cdot b\right)  \sim\Gamma\left(  1\cdot b\right)  $. In other words,
$\Gamma\left(  a\right)  \sim\Gamma\left(  b\right)  $ (since $a=\gamma\cdot
b$ and $b=1\cdot b$). Hence, Claim 2 is proved.
\end{proof}

Now, the map%
\begin{align*}
\overline{\Delta}:\mathcal{G}^{n}/\left(  C_{n}\times\mathcal{G}\right)   &
\rightarrow\mathcal{G}_{1}^{n}/C_{n},\\
\left(  C_{n}\times\mathcal{G}\text{-orbit of }a\right)   &  \mapsto\left(
C_{n}\text{-orbit of }\Delta\left(  a\right)  \right)
\end{align*}
is well-defined (by Claim 1), and so is the map%
\begin{align*}
\overline{\Gamma}:\mathcal{G}_{1}^{n}/C_{n}  &  \rightarrow\mathcal{G}%
^{n}/\left(  C_{n}\times\mathcal{G}\right)  ,\\
\left(  C_{n}\text{-orbit of }a\right)   &  \mapsto\left(  C_{n}%
\times\mathcal{G}\text{-orbit of }\Gamma\left(  a\right)  \right)
\end{align*}
(by Claim 2). We shall now prove that these two maps $\overline{\Delta}$ and
$\overline{\Gamma}$ are mutually inverse. We do this through two claims:

\begin{statement}
\textit{Claim 3:} We have $\overline{\Gamma}\circ\overline{\Delta
}=\operatorname*{id}$.
\end{statement}

\begin{proof}
[Proof of Claim 3.]Let $x\in\mathcal{G}^{n}/\left(  C_{n}\times\mathcal{G}%
\right)  $ be arbitrary. Thus, $x$ is the $C_{n}\times\mathcal{G}$-orbit of
some $a\in\mathcal{G}^{n}$. Consider this $a$. Thus, the definition of
$\overline{\Delta}$ shows that $\overline{\Delta}\left(  x\right)  $ is the
$C_{n}$-orbit of $\Delta\left(  a\right)  $. Hence, the definition of
$\overline{\Gamma}$ shows that $\overline{\Gamma}\left(  \overline{\Delta
}\left(  x\right)  \right)  $ is the $C_{n}\times\mathcal{G}$-orbit of
$\Gamma\left(  \Delta\left(  a\right)  \right)  $. We shall now show that
$\Gamma\left(  \Delta\left(  a\right)  \right)  $ and $a$ belong to the same
$C_{n}\times\mathcal{G}$-orbit. More specifically, we shall show that
$\Gamma\left(  \Delta\left(  a\right)  \right)  =\left(  g^{0},a_{\overline
{0}}^{-1}\right)  \cdot a$, where the pair $\left(  g^{0},a_{\overline{0}%
}^{-1}\right)  \in C_{n}\times\mathcal{G}$ acts on $a$ via the $C_{n}%
\times\mathcal{G}$-action on $\mathcal{G}^{n}$.

For each $k\in\left\{  1,2,\ldots,n\right\}  $, we have%
\begin{align*}
\left(  \Gamma\left(  \Delta\left(  a\right)  \right)  \right)  _{\overline
{k}} &  =\underbrace{\left(  \Delta\left(  a\right)  \right)  _{\overline{1}}%
}_{\substack{=a_{\overline{0}}^{-1}a_{\overline{1}}\\\text{(by
(\ref{pf.thm.hg-main.Del2}))}}}\cdot\underbrace{\left(  \Delta\left(
a\right)  \right)  _{\overline{2}}}_{\substack{=a_{\overline{1}}%
^{-1}a_{\overline{2}}\\\text{(by (\ref{pf.thm.hg-main.Del2}))}}}\cdot
\cdots\cdot\underbrace{\left(  \Delta\left(  a\right)  \right)  _{\overline
{k}}}_{\substack{=a_{\overline{k-1}}^{-1}a_{\overline{k}}\\\text{(by
(\ref{pf.thm.hg-main.Del2}))}}}\\
&\qquad\qquad\left(
\text{by (\ref{pf.thm.hg-main.Gam2}), applied to }\Delta\left(  a\right)
\text{ and }k
\text{ instead of }a\text{ and }j
\right)  \\
&  =a_{\overline{0}}^{-1}a_{\overline{1}}\cdot a_{\overline{1}}^{-1}%
a_{\overline{2}}\cdot\cdots\cdot a_{\overline{k-1}}^{-1}a_{\overline{k}%
} \\
&=a_{\overline{0}}^{-1}a_{\overline{k}}\qquad\left(  \text{by the telescope
principle}\right)
\end{align*}
and%
\begin{align*}
\left(  \left(  g^{0},a_{\overline{0}}^{-1}\right)  \cdot a\right)
_{\overline{k}} &  =a_{\overline{0}}^{-1}a_{\overline{k-0}}\qquad\left(
\text{by (\ref{pf.thm.hg-main.actCG2}), applied to }i=0\text{ and
}h=a_{\overline{0}}^{-1}\text{ and }j=k\right)  \\
&  =a_{\overline{0}}^{-1}a_{\overline{k}}\qquad\left(  \text{since
}k-0=k\right)  .
\end{align*}
Comparing these two equalities, we see that%
\[
\left(  \Gamma\left(  \Delta\left(  a\right)  \right)  \right)  _{\overline
{k}}=\left(  \left(  g^{0},a_{\overline{0}}^{-1}\right)  \cdot a\right)
_{\overline{k}}\qquad\text{for each }k\in\left\{  1,2,\ldots,n\right\}  .
\]
In other words, each entry of the $n$-tuple $\Gamma\left(  \Delta\left(
a\right)  \right)  $ equals the corresponding entry of $\left(  g^{0}%
,a_{\overline{0}}^{-1}\right)  \cdot a$. Hence, $\Gamma\left(  \Delta\left(
a\right)  \right)  =\left(  g^{0},a_{\overline{0}}^{-1}\right)  \cdot a$.
Therefore, $\Gamma\left(  \Delta\left(  a\right)  \right)  $ and $a$ belong to
the same $C_{n}\times\mathcal{G}$-orbit. In other words, the $C_{n}%
\times\mathcal{G}$-orbit of $\Gamma\left(  \Delta\left(  a\right)  \right)  $
is precisely the $C_{n}\times\mathcal{G}$-orbit of $a$. In other words,
$\overline{\Gamma}\left(  \overline{\Delta}\left(  x\right)  \right)  =x$
(since $\overline{\Gamma}\left(  \overline{\Delta}\left(  x\right)  \right)  $
is the $C_{n}\times\mathcal{G}$-orbit of $\Gamma\left(  \Delta\left(
a\right)  \right)  $, whereas $x$ is the $C_{n}\times\mathcal{G}$-orbit of $a$).

Forget that we fixed $x$. We thus have shown that $\overline{\Gamma}\left(
\overline{\Delta}\left(  x\right)  \right)  =x$ for each $x\in\mathcal{G}%
^{n}/\left(  C_{n}\times\mathcal{G}\right)  $. In other words, $\overline
{\Gamma}\circ\overline{\Delta}=\operatorname*{id}$. This proves Claim 3.
\end{proof}

\begin{statement}
\textit{Claim 4:} We have $\overline{\Delta}\circ\overline{\Gamma
}=\operatorname*{id}$.
\end{statement}

\begin{proof}
[Proof of Claim 4.]Let $x\in\mathcal{G}_{1}^{n}/C_{n}$ be arbitrary. Thus, $x$
is the $C_{n}$-orbit of some $a\in\mathcal{G}_{1}^{n}$. Consider this $a$.
Hence, the definition of $\overline{\Gamma}$ shows that $\overline{\Gamma
}\left(  x\right)  $ is the $C_{n}\times\mathcal{G}$-orbit of $\Gamma\left(
a\right)  $. Therefore, the definition of $\overline{\Delta}$ shows that
$\overline{\Delta}\left(  \overline{\Gamma}\left(  x\right)  \right)  $ is the
$C_{n}$-orbit of $\Delta\left(  \Gamma\left(  a\right)  \right)  $. We shall
now show that $\Delta\left(  \Gamma\left(  a\right)  \right)  =a$.

Let $k\in\left\{  1,2,\ldots,n\right\}  $. Then,
(\ref{pf.thm.hg-main.Del2}) (applied to $\Gamma\left(  a\right)$
and $k$ instead of $a$ and $j$) yields
\begin{align*}
\left(  \Delta\left(  \Gamma\left(  a\right)  \right)  \right)  _{\overline
{k}} &  =\left(  \Gamma\left(  a\right)  \right)  _{\overline{k-1}}^{-1}%
\cdot\left(  \Gamma\left(  a\right)  \right)  _{\overline{k}} \\
&  =\left(  a_{\overline{1}}a_{\overline{2}}\cdots a_{\overline{k-1}}\right)
^{-1}\cdot\left(  a_{\overline{1}}a_{\overline{2}}\cdots a_{\overline{k}%
}\right)
\end{align*}
(since (\ref{pf.thm.hg-main.Gam2a}) shows that
$\left(  \Gamma\left(
a\right)  \right)  _{\overline{k}}=a_{\overline{1}}a_{\overline{2}}\cdots
a_{\overline{k}}$ and
$\left(  \Gamma\left(  a\right)  \right)  _{\overline{k-1}%
}=a_{\overline{1}}a_{\overline{2}}\cdots a_{\overline{k-1}}$
(because $k-1\in\left\{  0,1,\ldots,n\right\}$)).
Thus,
\begin{align*}
\left(  \Delta\left(  \Gamma\left(  a\right)  \right)  \right)  _{\overline
{k}} &  =
\left(  a_{\overline{1}}a_{\overline{2}}\cdots a_{\overline{k-1}}\right)
^{-1}\cdot\left(  a_{\overline{1}}a_{\overline{2}}\cdots a_{\overline{k}%
}\right) \\
&  =a_{\overline{k}}\qquad\left(  \text{since }a_{\overline{1}}a_{\overline
{2}}\cdots a_{\overline{k}}=\left(  a_{\overline{1}}a_{\overline{2}}\cdots
a_{\overline{k-1}}\right)  \cdot a_{\overline{k}}\right)  .
\end{align*}

Forget that we fixed $k$.
We have thus shown that
$\left(  \Delta\left(  \Gamma\left(  a\right)  \right)  \right)  _{\overline
{k}} = a_{\overline{k}}$ for each $k \in \set{1,2,\ldots,n}$.
In other words, each entry of the $n$-tuple $\Delta\left(  \Gamma\left(
a\right)  \right)  $ equals the corresponding entry of $a$. Hence,
$\Delta\left(  \Gamma\left(  a\right)  \right)  =a$. Therefore, $\overline
{\Delta}\left(  \overline{\Gamma}\left(  x\right)  \right)  =x$ (since
$\overline{\Delta}\left(  \overline{\Gamma}\left(  x\right)  \right)  $ is the
$C_{n}$-orbit of $\Delta\left(  \Gamma\left(  a\right)  \right)  $, whereas
$x$ is the $C_{n}$-orbit of $a$).

Forget that we fixed $x$. We thus have shown that $\overline{\Delta}\left(
\overline{\Gamma}\left(  x\right)  \right)  =x$ for each $x\in\mathcal{G}%
_{1}^{n}/C_{n}$. In other words, $\overline{\Delta}\circ\overline{\Gamma
}=\operatorname*{id}$. This proves Claim 4.
\end{proof}

Combining Claim 3 with Claim 4, we see that the maps $\overline{\Delta}$ and
$\overline{\Gamma}$ are mutually inverse. In particular, $\overline{\Delta
}:\mathcal{G}^{n}/\left(  C_{n}\times\mathcal{G}\right)  \rightarrow
\mathcal{G}_{1}^{n}/C_{n}$ is an invertible map, i.e., a bijection. Hence, by
the bijection principle, we have%
\[
\left\vert \mathcal{G}^{n}/\left(  C_{n}\times\mathcal{G}\right)  \right\vert
=\left\vert \mathcal{G}_{1}^{n}/C_{n}\right\vert =\frac{1}{n}\sum_{d\mid
n}\phi\left(  \frac{n}{d}\right)  \left[  \mathcal{G}\div\frac{n}{d}\right]
\abs{ \mathcal{G} }  ^{d-1}
\]
(by Theorem \ref{thm.main}). This proves Theorem \ref{thm.hg-main}.
\end{proof}

\section{\label{sect.Kprod}$K$-product necklaces}

We can generalize our earlier definition of identity-product necklaces
by replacing the identity element with an arbitrary conjugacy-invariant
subset.

For this entire section, we fix a subset $K$ of $\G$ that is closed under conjugation (i.e., satisfies $gkg^{-1} \in K$ for all $g \in \G$ and $k \in K$).
Equivalently, $K$ shall be a union of conjugacy classes of $\G$.
Examples of such subsets $K$ are single conjugacy classes, normal subgroups, complements of normal subgroups, singleton subsets $\set{h}$ where $h$ lies in the center of $\G$, and the whole group $\G$.

\begin{definition}
    For any nonnegative integer $n$, we let $\G_K^n$ be the set of all \emph{$K$-product $n$-tuples}, defined by
    \begin{align*}
        \G_K^n &:= \set{a = (a_1, a_2, \dots, a_n) \in \G^n : a_1 a_2 a_3 \cdots a_n \in K}.
    \end{align*}
\end{definition}
The cyclic group $C_n$ acts on $\G^n$ by cyclic rotation, as we know.
We claim that the subset $\G_K^n$ is preserved under this action also.
That is, for any $K$-product $n$-tuple $(a_1, a_2, \dots, a_n)$, its cyclic rotation $(a_n, a_1, a_2, \dots, a_{n-1})$ is again a $K$-product $n$-tuple, since its product
\begin{align*}
    a_n a_1 a_2 \cdots a_{n-1}
    &= a_n \underbrace{a_1 a_2 \cdots a_{n-1} a_n}_{\in K} a_n^{-1} \in a_n K a_n^{-1} \subseteq K
\end{align*}
(since $K$ is closed under conjugation).
Thus $\G_K^n$ is a $C_n$-set. In a similar fashion as with identity-product necklaces, orbits of this action of $C_n$ on $\G_K^n$ will be called \emph{$K$-product $n$-necklaces}.

\begin{definition}
    For any integer $n$, let us define the number $[\G \div K \div n]$ to be
    \begin{align*}
        [\G \div K \div n] &:= \abs{\set{a \in \G : a^n \in K}}.
    \end{align*}
\end{definition}

This generalizes $[\G \div n]$, since $[\G \div n] = [\G \div \set{1} \div n]$.

With this generalization to $K$-product $n$-necklaces, we get an analogue to Theorem~\ref{thm.main}:
\begin{theorem}
    \label{thm.Kgeneralized}
    Let $n$ be a positive integer. The number $\abs{\G_K^n / C_n}$ of $K$-product $n$-necklaces is
    \begin{align*}
        \abs{\G_K^n / C_n} &= \frac{1}{n}\sum_{d \mid n} \phi\left(\frac{n}{d}\right)\left[\G \div K \div \frac{n}{d}\right] \abs{\G}^{d-1} .
    \end{align*}
\end{theorem}

The proof of this analogue of Theorem~\ref{thm.main} is the same except for the following changes: replace all  $\G_1$ by $\G_K$ and replace all terms of the form $[\G \div m]$ by $[\G \div K \div m]$.

For example, the analogous versions of the definitions of $P_k$ and $A_k$ as well as Lemma~\ref{lem.sizeOfAi} are, respectively
\begin{align}
    P_k &:= \set{a \in \G_K^n : g^k \cdot a = a, \text{ but } g^l \cdot a \neq a \text{ for all } 0 < l < k}. \\
    A_k &:= \set{a \in \G_K^n : g^k \cdot a = a}.
\end{align}
\begin{lemma}
    For any $k \mid n$, we have
    \begin{align}
        \abs{A_k} = \left[\G \div K \div \frac{n}{k}\right] \abs{\G}^{k-1} .
    \end{align}
\end{lemma}

This generalization also subsumes the classical formula \eqref{eq.neck}. Indeed, if we take $\G = C_q$ and $K = \G$, then Theorem~\ref{thm.Kgeneralized} recovers the count of $n$-necklaces over a size-$q$ set.

With the same replacement of all mentions of $\G_1$ and $[\G \div n]$ with their $K$-product analogues, we can also obtain an analogous $K$-product version of Theorem~\ref{thm.aperiodic} which counts the number of aperiodic $K$-product $n$-necklaces:
\begin{theorem}
    Let $n$ be a positive integer. The number $\frac{\abs{P_n}}{n}$ of aperiodic $K$-product $n$-necklaces is
    \begin{align*}
        \frac{\abs{P_n}}{n} &= \frac{1}{n}\sum_{d \mid n} \mu\left(\frac{n}{d}\right) \left[\G \div K \div \frac{n}{d}\right] \abs{\G}^{d-1} .
    \end{align*}
\end{theorem}

We leave it as an open question as to whether or not the bijection between identity-product necklaces and homogeneous necklaces has some analogue in the $K$-product setting as well.

\section{\label{sect.irpol}Irreducible polynomials}

As we briefly mentioned in the introduction, aperiodic necklaces over a size-$q$ set are related to irreducible degree-$n$ polynomials over the finite field $\mathbb{F}_q$ when $q$ is a prime power.
This relation can be restricted to aperiodic identity-product $n$-necklaces (over the additive group of $\mathbb{F}_q$, so the identity is actually $0$) on the one hand, and to the irreducible degree-$n$ polynomials with second-highest coefficient $0$ on the other.
We shall now elaborate on this.

\subsection{Reminders on finite fields and their Galois groups}

Let $q$ be a prime power.
Let $K / F$ be a degree-$n$ field extension, where $F$ is the finite field with $q$ elements.
Thus, $K$ is a finite field with $q^n$ elements.
Note that $K / F$ is a Galois extension (like any extension of finite fields).

It is well-known that the map
\begin{align*}
    \Frob : K &\to K, \\
    x &\mapsto x^q
\end{align*}
is an $F$-algebra automorphism of $K$.
This is called the \emph{Frobenius automorphism over $\mathbb{F}_q$}, and generates the Galois group of the field extension $K/F$.
Its $k$-th power $\Frob^k$ sends each $x \in K$ to $x^{q^k}$.
In particular, its $n$-th power $\Frob^n$ is the identity map, since each $x \in K$ satisfies $x^{q^n} = x^{\abs{K}} = x$.
Thus, the Galois group $\Gal(K / F)$ consists of the $n$ distinct automorphisms $\Frob^0, \Frob^1, \ldots, \Frob^{n-1}$.

A significant role in Galois theory is played by the trace:

\begin{definition}
    The \emph{trace map} is the $F$-linear map $\tr : K \to K$ defined by
    \begin{equation}
        \tr := \Frob^0 + \Frob^1 + \cdots + \Frob^{n-1}.
        \label{eq.def.trace.0}
    \end{equation}

    The \emph{trace} of an element $x \in K$ is defined to be the image $\tr(x)$ of $x$ under this map. Explicitly,
    \begin{align}
    \tr(x) &= x + \Frob x + \Frob^2 x + \cdots + \Frob^{n-1} x
    \label{eq.def.trace.1} \\
    &= x + x^q + x^{q^2} + \cdots + x^{q^{n-1}} .
    \label{eq.def.trace.2}
\end{align}
\end{definition}

Notably, $\tr \circ \, \Frob = \tr$, since right-multiplication by $\Frob$ merely cyclically permutes the addends on the right hand side of \eqref{eq.def.trace.0} (since $\Frob^n = \id = \Frob^0$) but does not affect the sum as a whole.
Similarly, $\Frob \circ \tr = \tr$.
The latter equality shows that each $x \in K$ satisfies $\Frob\tup{\tr\tup{x}} = \tr\tup{x}$, and thus $\tr\tup{x}$ is fixed by the Galois group $\Gal(K/F)$ (since all elements of $\Gal(K/F)$ are powers of $\Frob$).
Since $K / F$ is a Galois extension, this shows that
\begin{equation}
\tr(x) \in F \qquad \text{ for all } x \in K.
\label{eq.def.trace.in-F}
\end{equation}

\begin{definition}
    We say that two elements of $K$ are \emph{conjugate} if they have the same minimal polynomial over $F$. Equivalently, two elements $x,y \in K$ are conjugate if $x = A\tup{y}$ for some $A \in \Gal(K/F)$.
\end{definition}

If $x$ and $y$ are two conjugate elements of $K$, then their traces are equal: $\tr(x) = \tr(y)$.
Indeed, from $\tr \circ \, \Frob = \tr$, we obtain
\begin{align}
    \tr \circ \, \Frob^k = \tr \qquad \text{for all } k \in \mathbb{Z}.
    \label{eq.trace.trFk}
\end{align}
But this means that $\tr \circ \, A = \tr$ for each $A \in \Gal(K/F)$ (since all elements of $\Gal(K/F)$ are powers of $\Frob$).
Hence, if $x = A\tup{y}$ for some $A \in \Gal(K/F)$, then $\tr\tup{x} = \tr\tup{A\tup{y}} = \tup{\tr \circ\, A} \tup{y} = \tr\tup{y}$.

Hence, we define the trace of an entire conjugacy class to be the trace of any of its elements.

\subsection{Normal bases and the bijection $\Phi$}

A \emph{normal basis} of $K$ (over $F$) means a basis of $K$ as an $F$-vector space that has the form $\tup{\theta^{q^0}, \theta^{q^1}, \ldots, \theta^{q^{n-1}}}$, where $\theta \in K$.
A classical fact (the normal basis theorem for finite fields) says:

\begin{theorem}
\label{thm.nbf}
There exists a normal basis of $K$ over $F$.
\end{theorem}

Proofs of this theorem can be found in \cite[Theorem 2.35]{LidNie}, \cite[Appendix B]{Murphy-FF}, \cite[Theorem 1.15]{Hou}, \cite[Corollary 2.4.6]{Gao}
and various other texts on finite fields.
It is in fact a particular case of the normal basis theorem for arbitrary Galois extensions, which appears (e.g.) in \cite[Theorem 5.18]{Milne-FT}, \cite{Blessenohl} etc.






Henceforth, we fix a normal basis $\tup{\theta^{q^0}, \theta^{q^1}, \ldots, \theta^{q^{n-1}}}$ of $K$ over $F$.

Following \cite[Section 7.6.2]{Reutenauer}, we introduce the map
\begin{align*}
    \Phi : F^n &\to K, \\
    \tup{a_0, a_1, \ldots, a_{n-1}} &\mapsto
    a_0 \theta^{q^0} + a_1 \theta^{q^1} + \cdots + a_{n-1} \theta^{q^{n-1}}.
\end{align*}
This is a bijection, since $\tup{\theta^{q^0}, \theta ^{q^1}, \ldots, \theta^{q^{n-1}}}$ is a basis of $K$ over $F$.
Note that $\Frob^n\tup{\theta} = \theta$, since $\Frob^n = \id$.

Recall that the cyclic group $C_n$ acts on the set $F^n$ by cyclically rotating $n$-tuples.
We also let $C_n$ act on the set $K$ by having the generator $g$ act as the Frobenius automorphism $\Frob$.
This is a well-defined group action, since $\Frob^n = \id$. Now we claim the following:

\begin{lemma}
    \label{lem.phi-equivariant}
    The bijection $\Phi$ is an isomorphism of $C_n$-sets.
\end{lemma}

\begin{proof}
    We must prove that $\Phi$ is $C_n$-equivariant (since we already know $\Phi$ to be a bijection).
    Thus, we need to show that $\Phi\left(h \cdot a\right) = h \cdot \Phi\tup{a}$ for any $a \in F^n$ and any $h \in C_n$.
    Since the group $C_n$ is generated by $g$, it suffices to show this for $h = g$.
    So we must prove that $\Phi\left(g \cdot a\right) = g \cdot \Phi\tup{a}$ for any $a \in F^n$.

    Consider the standard basis $\tup{e_0, e_1, \ldots, e_{n-1}}$ of the $F$-vector space $F^n$,
    where $e_i$ is the $n$-tuple $\tup{0,0,\dots,0,1,0,0,\dots,0}$ with $i$ many $0$'s preceding the entry $1$.
    We also set $e_n := e_0$.
    The $C_n$-action on $F^n$ then satisfies $g \cdot e_i = e_{i+1}$ for each $i \in \set{0,1,\ldots,n-1}$.

    We must prove that $\Phi\left(g \cdot a\right) = g \cdot \Phi\tup{a}$ for any $a \in F^n$.
    Since $\Phi$ is $F$-linear, it suffices to prove this for all basis vectors $a = e_i$ with $i \in \set{0,1,\ldots,n-1}$.
    That is, we must show that $\Phi\left(g \cdot e_i\right) = g \cdot \Phi\tup{e_i}$ for all $i \in \set{0,1,\ldots,n-1}$.

    Fix $i \in \set{0,1,\ldots,n-1}$.
    Then, $g \cdot e_i = e_{i+1}$ as we saw above.
    Hence, $\Phi\tup{g \cdot e_i} = \Phi\tup{e_{i+1}}$.
    But the definition of $\Phi$ yields $\Phi\tup{e_i} = \theta^{q^i} = \Frob^i\tup{\theta}$, so that
    \begin{align}
        g \cdot \Phi\tup{e_i} = g \cdot \Frob^i\tup{\theta} = \Frob\tup{\Frob^i\tup{\theta}} = \Frob^{i+1}\tup{\theta}.
        \label{pf.lem.phi-equivariant.gPei}
    \end{align}
    Moreover, we claim that $\Phi\tup{e_{i+1}} = \Frob^{i+1}\tup{\theta}$. Indeed, if $i+1 < n$, then this follows from the definition of $\Phi$ (similarly to $\Phi\tup{e_i} = \Frob^i\tup{\theta}$), whereas in the case $i+1 = n$, it follows from $\Phi\tup{e_n} = \Phi\tup{e_0} = \theta^{q^0} = \theta = \Frob^n\tup{\theta}$.
    Thus, $\Phi\left(g \cdot e_i\right) = \Phi\tup{e_{i+1}} = \Frob^{i+1}\tup{\theta}$.
    Comparing this with \eqref{pf.lem.phi-equivariant.gPei}, we obtain
    $\Phi\left(g \cdot e_i\right) = g \cdot \Phi\tup{e_i}$ as desired.
\end{proof}

\begin{corollary}
    \label{cor.phimap}
    Two $n$-tuples in $F^n$ lie in the same $C_n$-orbit if and only if their images under $\Phi$ are conjugate elements.
\end{corollary}

\begin{proof}
    Let $x$ and $y$ be two $n$-tuples in $F^n$.
    By Lemma~\ref{lem.phi-equivariant}, these tuples $x$ and $y$ lie in the same $C_n$-orbit if and only if their images $\Phi\tup{x}$ and $\Phi\tup{y}$ lie in the same $C_n$-orbit.
    But the $C_n$-orbits on $K$ are just the orbits of the Galois group $\Gal(K / F)$, since $C_n$ acts on $K$ by having $g$ act as the generator $\Frob$ of $\Gal(K / F)$.
    So $x$ and $y$ lie in the same $C_n$-orbit if and only if their images $\Phi\tup{x}$ and $\Phi\tup{y}$ lie in the same orbit of $\Gal(K / F)$, that is, are conjugate.
\end{proof}

Thus, $\Phi$ gives rise to a bijection from $F^n / C_n$ (that is, $n$-necklaces over $F$) to the set of conjugacy classes of elements of $K$ over $F$.
We denote this bijection by $\hat{\Phi}$.

\subsection{Zero-sum necklaces and zero-trace conjugacy classes}

We say that an $n$-tuple $a = \tup{a_1, a_2, \ldots, a_n} \in F^n$ is \emph{zero-sum} if $a_1 + a_2 + \cdots + a_n = 0$ in $F$.
The zero-sum $n$-tuples are just the identity-product $n$-tuples for the abelian group $\tup{F,+,0}$. We denote the set of those zero-sum $n$-tuples as $F_0^n$.

We say that an $n$-necklace $\mathbf{a} \in F^n / C_n$ is \emph{zero-sum} if its representatives $a = \tup{a_1, a_2, \ldots, a_n} \in \mathbf{a}$ are zero-sum $n$-tuples (i.e., satisfy $a_1 + a_2 + \cdots + a_n = 0$ in $F$).
This condition clearly does not depend on the choice of the representative, since cyclic rotation of an $n$-tuple does not change the sum of its entries.
The zero-sum $n$-necklaces are just the identity-product $n$-necklaces for the abelian group $\tup{F,+,0}$. Therefore, their set is $F_0^n / C_n$.

\begin{theorem}
    \label{thm.zero-sumIsTrace}
    A necklace $\mathbf{a} \in F^n / C_n$ is zero-sum if and only if the corresponding conjugacy class $\hat{\Phi}\tup{\mathbf{a}}$ has trace $0$.
\end{theorem}

\begin{proof}
    [Proof of Theorem \ref{thm.zero-sumIsTrace}]
    Let us first note that \eqref{eq.def.trace.2} yields $\tr(\theta) = \theta + \theta^q + \theta^{q^2} + \cdots + \theta^{q^{n-1}} \neq 0$, since $\tup{\theta, \theta^q, \theta^{q^2}, \ldots, \theta^{q^{n-1}}} = \tup{\theta^{q^0}, \theta ^{q^1}, \ldots, \theta^{q^{n-1}}}$ is a basis of $K$ over $F$.

    Fix a necklace $\mathbf{a}$, and let $a = \tup{a_0, a_1, \ldots, a_{n-1}}$ be a representative (i.e., one of the tuples in $\mathbf{a}$). Thus, $\hat{\Phi}\tup{\mathbf{a}}$ is the conjugacy class of
    \begin{align*}
    \Phi\tup{a}
    &= \Phi\tup{a_0, a_1, \ldots, a_{n-1}}
    = a_0 \theta^{q^0} + a_1 \theta^{q^1} + \cdots + a_{n-1} \theta^{q^{n-1}} \\
    &= \sum_{k=0}^{n-1} a_k \theta^{q^k}
    = \sum_{k=0}^{n-1} a_k \Frob^k\tup{\theta}.
    \end{align*}
    Hence, we can compute the trace of $\hat{\Phi}\tup{\mathbf{a}}$ as follows:
    \begin{align*}
        \tr(\hat{\Phi}\tup{\mathbf{a}}) &=
        \tr\left(\sum_{k=0}^{n-1} a_k \Frob^k\tup{\theta}\right) \\
        &= \sum_{k=0}^{n-1} a_k \tr \tup{\Frob^k(\theta)} \qquad \left(\text{since $\tr$ is $F$-linear}\right)\\
        &= \sum_{k=0}^{n-1} a_k \tup{\tr \circ\, \Frob^k}  (\theta) \\
        &= \sum_{k=0}^{n-1} a_k \tr (\theta) \qquad \left(\text{by \eqref{eq.trace.trFk}}\right) \\
        &= (a_0 + a_1 + \cdots + a_{n-1}) \tr(\theta).
    \end{align*}
    Thus, $\tr(\hat{\Phi}\tup{\mathbf{a}}) = 0$ if and only if $a_0 + a_1 + \cdots + a_{n-1} = 0$ (since $\tr(\theta) \neq 0$).
    In other words, $\hat{\Phi}\tup{\mathbf{a}}$ has trace $0$ if and only if $\mathbf{a}$ is a zero-sum necklace.
\end{proof}

\begin{corollary}
    The number of conjugacy classes of elements of $K$ over $F$ with trace $0$ is equal to $\abs{F^n_0 / C_n}$.
\end{corollary}
\begin{proof}
    Theorem~\ref{thm.zero-sumIsTrace} shows that the bijection $\hat{\Phi}$ restricts to a bijection from the set $F^n_0 / C_n$ of zero-sum $n$-necklaces to the set of conjugacy classes of elements of $K$ over $F$ with trace $0$.
    Thus, the two sets are equinumerous.
\end{proof}

Thus, we can use Theorem~\ref{thm.main} to count the conjugacy classes with trace $0$.

We shall next restrict ourselves to zero-sum $n$-tuples and $n$-necklaces with certain periodicity properties. We recall the notions of ``aperiodic'' and ``smallest period'' we introduced in Section~\ref{sec.aperiodic}.
We define $P_n$ to be the set of all zero-sum $n$-tuples in $F^n_0$ whose smallest period is $n$.
This is a particular case of the definition of $P_k$ in \eqref{def.Pk}, applied to the group $\tup{F, +, 0}$ as $\G$.

\begin{proposition}
    \label{prop.conjClassesToNecklaces}
    The number of size-$n$ conjugacy classes of elements of $K$ with trace $0$ is equal to the number of aperiodic zero-sum $n$-necklaces.
\end{proposition}

\begin{proof}
    Since $\Phi: F^n \to K$ is a bijection, it preserves the size of any subset of $F^n$.
    Thus, the size of an $n$-necklace $\mathbf{a}$ as a set must be equal to the size of the conjugacy class $\hat{\Phi}(\mathbf{a})$. In particular, the $n$-necklaces with size $n$ as sets are bijectively mapped by $\hat{\Phi}$ to the size-$n$ conjugacy classes.
    But the former necklaces are precisely the aperiodic $n$-necklaces.
    Thus, by restricting $\hat{\Phi}$ to these necklaces, we obtain a bijection from the aperiodic $n$-necklaces to the size-$n$ conjugacy classes.

    Restricting this bijection further to the zero-sum aperiodic $n$-necklaces, we obtain precisely the size-$n$ conjugacy classes with trace $0$ as their images because of Theorem~\ref{thm.zero-sumIsTrace}.
\end{proof}

By Theorem~\ref{thm.aperiodic}, the number of aperiodic zero-sum $n$-necklaces is $\frac{\abs{P_n}}{n}$.

\begin{corollary}
    There are $\abs{P_n}$ elements in $K$ with trace $0$ that belong to size-$n$ conjugacy classes.
\end{corollary}

\begin{proof}
    By Proposition~\ref{prop.conjClassesToNecklaces}, the number of size-$n$ conjugacy classes of elements with trace $0$ in $K$ is the number of aperiodic zero-sum $n$-necklaces; but this number is $\frac{\abs{P_n}}{n}$ by Theorem~\ref{thm.aperiodic}.
    Since these classes are disjoint and have size $n$ each, they contain a total of $\frac{\abs{P_n}}{n} \cdot n = \abs{P_n}$ many elements.
\end{proof}

\subsection{The viewpoint of irreducible polynomials}

Proposition~\ref{prop.conjClassesToNecklaces} leads to an enumerative result about irreducible polynomials over $F$. To obtain it, we need a lemma that connects the trace of an element of $K$ with its minimal polynomial:

\begin{lemma}
    \label{lem.traceIsCoefficient}
    Let $z \in K$ be an element whose minimal polynomial over $F$ has degree $n$.
    Then, the second-highest coefficient of this polynomial is $-\tr(z)$.
\end{lemma}

This is a special case (where the minimal polynomial is of degree $n$) of \cite[Corollary 5.45]{Milne-FT} or of \cite[Section 14.2, Exercise 18 (d)]{Dummit-Foote}. 

\begin{corollary}
    \label{cor.num-irrpols}
    The number of monic irreducible polynomials of degree $n$ over $F$ with second-highest coefficient equal to $0$ is equal to the number of aperiodic zero-sum $n$-necklaces over $F$.
\end{corollary}

\begin{proof}
    Starting from Proposition~\ref{prop.conjClassesToNecklaces}, we know that the number of aperiodic zero-sum $n$-necklaces over $F$ is equal to the number of size-$n$ conjugacy classes with trace $0$.

    Each monic irreducible polynomial $f$ of degree $n$ over $F$ has $n$ distinct zeroes in $K$ (since $F$ is a perfect field and $K$ is its only degree-$n$ extension up to isomorphism,
    and since $K/F$ is a Galois extension).
    These $n$ zeroes form a conjugacy class, and $f$ is their minimal polynomial.

    Thus, we can define a map
    \begin{align*}
        \sigma &: \set{\text{monic irreducible polynomials of degree $n$ over $F$}} \\
        &\to \set{\text{size-$n$ conjugacy classes on $K$}}
    \end{align*}
    that sends a polynomial to the conjugacy class consisting of its $n$ zeroes in $K$.
    The inverse $\sigma^{-1}$ of this map sends
    \begin{align*}
        \{z_1, z_2, \dots, z_n\} \mapsto \underbrace{(x - z_1)(x - z_2) \cdots (x - z_n)}_{\substack{\in F[x] \\ \text{(since this polynomial is preserved} \\ \text{by the Galois group } \Gal(K/F) \text{)}}}.
    \end{align*}
    Thus $\sigma$ is a bijection.
    
    Via $\sigma^{-1}$, each size-$n$ conjugacy class with trace $0$ corresponds to a monic irreducible degree-$n$ polynomial over $F$, and by Lemma~\ref{lem.traceIsCoefficient} that polynomial must have second-highest coefficient $0$; and this is a one-to-one correspondence.
    Thus the number of monic irreducible polynomials of degree $n$ over $F$ with second-highest coefficient equal to $0$
    is equal to the number of size-$n$ conjugacy classes on $K$ with trace $0$;
    but the latter number is in turn equal to the number of aperiodic zero-sum $n$-necklaces over $F$.
\end{proof}

Combining Corollary~\ref{cor.num-irrpols} with the explicit formula in Theorem~\ref{thm.aperiodic} for the number of aperiodic zero-sum $n$-necklaces over $F$, we obtain a formula for the number of monic irreducible polynomials of degree $n$ over $F$ with second-highest coefficient equal to $0$;
this recovers \cite[Corollary 2.8]{Yucas}.

\section{\label{sect.oeis}Sequences}

The numbers $\abs{\G_1^n / C_n}$ computed in Theorem~\ref{thm.main} give rise to several integer sequences, some of which are present in the OEIS (On-Line Encyclopedia of Integer Sequences, \cite{oeis}).
The simplest way to obtain a sequence is to fix the group $\G$ and to vary $n$.

Setting the group $\G$ to be the cyclic group $C_2$, and varying the length $n$ of the necklaces being counted, we obtain the OEIS sequence \seqnum{A000013}, which counts the equivalence classes of $n$-bead binary necklaces with beads of $2$ colors with respect to swapping all colors (the equivalence follows from Theorem~\ref{thm.hg-main}):
\[
    \abs{(C_2)_1^n / C_n} = \frac{1}{n}\sum_{d \mid n} \phi\left(\frac{n}{d}\right) \gcd\left(2, \frac{n}{d}\right)2^{d-1} .
\]
See Table~\ref{tab.1} for the first few values.
\begin{table}[H]
\centering
\begin{tabular}[c]{ |c||c|c|c|c|c|c|c|c|c| } 
 \hline
 $\vphantom{\int_g^b}$ $n$                  & 1 & 2 & 3 & 4 & 5 & 6 & 7 & 8 & 9 \\ \hline 
 $\vphantom{\int_g^b}$ $\abs{(C_2)_1^n / C_n}$ & 1 & 2 & 2 & 4 & 4 & 8 & 10 & 20 & 30
 \\ \hline
\end{tabular}
\caption{The numbers $\abs{(C_2)_1^n / C_n}$ for $n \leq 9$.}
\label{tab.1}
\end{table}
These numbers $\abs{(C_2)_1^n / C_n}$ appear under the name of $S\tup{n-1}$ in \cite{Sloane}.

We can also vary $n$ and $\G$ in lockstep.
For instance, we can let $\G$ be the cyclic group $C_n$ and count identity-product $n$-necklaces.
This gives us the OEIS sequence \seqnum{A130293}, which counts the number of necklaces of $n$ colors up to cyclic shifting of the colors (which are labelled $1, 2, \ldots, n$) \textbf{and} cyclic permutation of beads (again by Theorem~\ref{thm.hg-main}):
\begin{align*}
    \abs{(C_n)_1^n / C_n}
    &= \frac{1}{n}\sum_{d \mid n} \phi\left(\frac{n}{d}\right) \gcd\left(n, \frac{n}{d}\right) n^{d-1} \\
    &= \frac{1}{n}\sum_{d \mid n} \phi\left(\frac{n}{d}\right) \left(\frac{n}{d}\right) n^{d-1}
    \qquad \left(\text{since }\gcd\left(n, \frac{n}{d}\right) = \frac{n}{d} \right) \\
    &= \sum_{d \mid n} \phi\left(\frac{n}{d}\right) \frac{n^{d-1}}{d} .
\end{align*}
See Table~\ref{tab.2}.
\begin{table}[H]
\centering
\begin{tabular}[c]{ |c||c|c|c|c|c|c|c|c|c| } 
 \hline
 $\vphantom{\int_g^b}$ $n$                  & 1 & 2 & 3 & 4 & 5 & 6 & 7 & 8 & 9 \\ \hline 
 $\vphantom{\int_g^b}$ $\abs{(C_n)_1^n / C_n}$ & 1 & 2 & 5 & 20 & 129 & 1316 & 16813 & 262284 & 4783029
 \\ \hline
\end{tabular}
\caption{The numbers $\abs{(C_n)_1^n / C_n}$ for $n \leq 9$.}
\label{tab.2}
\end{table}

Alternatively, letting $\G$ be the cyclic group $C_{n+1}$ and again counting identity-product $n$-necklaces, we obtain the integer sequence \seqnum{A121774}. This counts the number of $n$-bead necklaces with $n+1$ colors divided by $n+1$:
\begin{align*}
    \abs{(C_{n+1})_1^n / C_n}
    &= \frac{1}{n}\sum_{d \mid n} \phi\left(\frac{n}{d}\right)
    \underbrace{\gcd\left(n+1, \frac{n}{d}\right)}_{\substack{= 1 \\ \text{(since } n+1 \text{ is coprime to $n$)}}}
    (n+1)^{d-1} \\
    &= \frac{1}{n}\sum_{d \mid n} \phi\left(\frac{n}{d}\right) (n+1)^{d-1}.
\end{align*}
See Table~\ref{tab.3}.
\begin{table}[H]
\centering
\begin{tabular}[c]{ |c||c|c|c|c|c|c|c|c|c| } 
 \hline
 $\vphantom{\int_g^b}$ $n$                  & 1 & 2 & 3 & 4 & 5 & 6 & 7 & 8 & 9 \\ \hline 
 $\vphantom{\int_g^b}$ $\abs{(C_{n+1})_1^n / C_n}$ & 1 & 2 & 6 & 33 & 260 & 2812 & 37450 & 597965 & 11111134
 \\ \hline
\end{tabular}
\caption{The numbers $\abs{(C_{n+1})_1^n / C_n}$ for $n \leq 9$.}
\label{tab.3}
\end{table}
A similar simplified formula can be obtained for $\abs{\G_1^n / C_n}$ whenever the size $\abs{\G}$ of $\G$ is coprime to $n$; but the result (which equals the number of $n$-bead necklaces with $\abs{\G}$ colors divided by $\abs{\G}$) can also be proved by a simple direct bijection, which we invite the reader to find.

\section{\label{sect.open}Further directions}

A variation on the set-up of the problem is to change the group acting on $\G_1^n$. A natural candidate for this is the dihedral group $D_n$, since $\G_1^n$ is already closed under the reflection $r$ defined by
\begin{align*}
    r \cdot (a_1, \dots, a_n) &:= (a_n^{-1}, \dots, a_1^{-1}) .
\end{align*}
Here, again, we can ask for the number of orbits, noticing that the analogous necklaces on non-group alphabets are known as \emph{bracelets} and have been counted long ago \cite{ZelZel}. 
When $\G$ is abelian, we can also have the symmetric group $S_n$ (or a subgroup thereof) act on $\G_1^n$ by permuting the entries.\footnote{Cf.\ also the notion of \emph{decimation classes} considered in \cite{Turner-etal} and \cite{Baczkowski-etal}.}

Generalizing in a different direction but in a similar vein, we can consider $n$-tuples of matrices whose product has a given trace:

\begin{definition}
    Let $R$ be a commutative ring.
    Consider $n$-tuples of $m \times m$-matrices over $R$.
    We define the set of \emph{$h$-trace-product $n$-tuples} $\prescript{h}{}{\mathcal{M}^n_m}$ to be 
    \begin{align*}
        \prescript{h}{}{\mathcal{M}^n_m}
        &:= \set{A = (A_1, A_2, \dots, A_n) : \tr(A_1 A_2 \cdots A_n) = h}.
    \end{align*}
\end{definition}
Since $\tr(AB) = \tr(BA)$, this set $\prescript{h}{}{\mathcal{M}^n_m}$ is a $C_n$-set and therefore we can attempt to count the number of orbits $\abs{\prescript{h}{}{\mathcal{M}^n_m} / C_n}$ when $R$ is finite.

\bibliographystyle{jis.bst}
\bibliography{Refs.bib}

\end{document}